\documentclass[11pt]{amsart}
\pdfoutput=1

  \pdfpageattr{/Group <</S /Transparency /I true /CS /DeviceRGB>>}

\usepackage[T1]{fontenc}
\usepackage[english]{babel}
\usepackage{hyperref}
\usepackage{geometry}
\usepackage{pdflscape}
\usepackage{amsfonts,amsmath,amssymb,amsthm}
\usepackage{tikz}
\usepackage{pgfplots}
\usepackage{caption}
\usepackage[labelfont=rm]{subcaption}
\usepackage{xspace}
\usepackage{booktabs,tabularx}
\usepackage{dcolumn}
\usepackage{listings}
\usepackage{marvosym}
\usepackage{microtype}

\usepackage{placeins}
\usepackage{longtable}

\microtypesetup{tracking,kerning,spacing}
\microtypecontext{spacing=nonfrench}

\usetikzlibrary{calc, 3d}
\newlength{\mywidth}
\newlength{\myheight}

\definecolor{links}{rgb}{.2,.1,.5}
\definecolor{cites}{rgb}{.5,.1,.2}
\hypersetup{colorlinks,linkcolor=links,citecolor=cites}

\newcolumntype{d}[1]{D{.}{.}{#1} }

\theoremstyle{plain}
\newtheorem{ruleofthumb}{Rule of Thumb}
\theoremstyle{definition}
\newtheorem{remark}{Remark}

\newcommand\scip{\texttt{SCIP}\xspace}
\newcommand\cplex{\texttt{CPLEX}\xspace}
\newcommand\gurobi{\texttt{Gurobi}\xspace}
\newcommand\miplib{\texttt{MIPLIB}\xspace}
\newcommand\polymake{\texttt{polymake}\xspace}
\newcommand\cdd{\texttt{cdd}\xspace}
\newcommand\ppl{\texttt{ppl}\xspace}
\newcommand\lrs{\texttt{lrs}\xspace}
\newcommand\bb{\texttt{bb}\xspace} % short for beneath_beyond
\newcommand\ftit{\texttt{4ti2}\xspace}
\newcommand\projection{\texttt{projection}\xspace}
\newcommand\latte{\texttt{LattE}\xspace}
\newcommand\porta{\texttt{porta}\xspace}
\newcommand\azove{\texttt{azove}\xspace}
\newcommand\lattem{\texttt{LattE macchiato}\xspace}
\newcommand\normaliz{\texttt{normaliz}\xspace}
\newcommand\proj{\texttt{proj}\xspace}
\newcommand\barvinok{\texttt{barvinok}\xspace}
\newcommand\cgal{\texttt{CGAL}\xspace}
\newcommand\nmz{\texttt{nmz}\xspace} % short for normaliz
\newcommand\gmp{\texttt{GMP}\xspace} % short for normaliz
\newcommand\bbox{\texttt{bbox}\xspace} % short for bounding_box
\newcommand\symcdd{\texttt{cdd$_{\text{\texttt{sym}}}$}\xspace}
\newcommand\symlrs{\texttt{lrs$_{\text{\texttt{sym}}}$}\xspace}
\newcommand\symppl{\texttt{ppl$_{\text{\texttt{sym}}}$}\xspace}
\newcommand\symbb{\texttt{bb$_{\text{\texttt{sym}}}$}\xspace}
\newcommand\sympol{\texttt{sympol}\xspace}
\newcommand\perl{\texttt{Perl}\xspace}
\newcommand\cplus{\texttt{C++}\xspace}

\definecolor{mycolor1}{rgb}{0.16863,0.50588,0.33725}%

\pgfplotsset{pplstyle/.style={color=orange,solid,mark=pentagon*,mark options={solid,fill=orange}}}
\pgfplotsset{bbstyle/.style={color=red,solid,mark=*,mark options={solid,fill=red}}}
\pgfplotsset{cddstyle/.style={color=blue,solid,mark=square*,mark options={solid,fill=blue}}}
\pgfplotsset{lrsstyle/.style={color=mycolor1,solid,mark=triangle*,mark options={solid,fill=mycolor1}}}
\pgfplotsset{nmzstyle/.style={color=black,solid,mark=diamond*,mark options={solid,fill=black}}}
\pgfplotsset{nmzstylealt/.style={color=black,solid,mark=diamond*,mark options={solid,fill=white}}}
\pgfplotsset{portastyle/.style={color=gray,solid,mark=10-pointed star,mark options={solid,fill=gray}}}
\pgfplotsset{azovestyle/.style={color=gray,solid,mark=10-pointed star,mark options={solid,fill=gray}}}
\pgfplotsset{bboxstyle/.style={color=red,solid,mark=*,mark options={solid,fill=red}}}
\pgfplotsset{projstyle/.style={color=orange,solid,mark=pentagon*,mark options={solid,fill=orange}}}
\pgfplotsset{ftitstyle/.style={color=mycolor1,solid,mark=triangle*,mark options={solid,fill=mycolor1}}}
\pgfplotsset{lattestyle/.style={color=blue,solid,mark=square*,mark options={solid,fill=blue}}}
\pgfplotsset{lattestylealt/.style={color=blue,solid,mark=square*,mark options={solid,fill=white}}}

\pgfplotsset{
legend image code/.code={
\draw[mark repeat=2,mark phase=2]
plot coordinates {
(0cm,0cm)
(0.15cm,0cm)        %% default is (0.3cm,0cm)
(0.3cm,0cm)         %% default is (0.6cm,0cm)
};%
}
}

\newcommand\RR{{\mathbb R}}
\newcommand\ZZ{{\mathbb Z}}

\newcommand\norm[1]{\left\lVert#1\right\rVert}
\newcommand\SetOf[2]{\left\{\left.#1\vphantom{#2}\ \right|\ #2\vphantom{#1}\right\}}
\newcommand\smallSetOf[2]{\{{#1}\,|\,{#2}\}}
\newcommand{\Cut}{\mathsf{Cut}}
\newcommand{\conv}{\mathsf{conv}}
\newcommand{\cone}{\mathsf{cone}}

\newcommand{\textoom}{\ToTop\xspace}
\newcommand{\oom}{\multicolumn{1}{c}{\textoom}}
\newcommand{\nc}{\multicolumn{1}{c}{-}}

\title{Computing convex hulls and counting integer points with \polymake}

\author[Assarf, Gawrilow, Herr, Joswig, Lorenz, Paffenholz, Rehn]%
{Benjamin Assarf \and Ewgenij Gawrilow \and Katrin Herr \and Michael Joswig \and Benjamin Lorenz \and Andreas Paffenholz \and Thomas Rehn}

\address[Benjamin Assarf, Michael Joswig, Benjamin Lorenz]{
Institut f{\"u}r Mathematik,
 TU Berlin,
 Str.\ des 17. Juni 136, 10623 Berlin, Germany
}
\email{\href{mailto:<lastname>@math.tu-berlin.de}{<lastname>@math.tu-berlin.de}}

\address[Ewgenij Gawrilow]{
}
\email{\href{mailto:egawrilow@gmail.com}{egawrilow@gmail.com}}

\address[Katrin Herr, Andreas Paffenholz]{
  Fachbereich Mathematik,
  TU Darmstadt,
  Dolivostr. 15,
  64293 Darmstadt, Germany
}
\email{\href{mailto:herr@mathematik.tu-darmstadt.de}{herr@mathematik.tu-darmstadt.de}, \href{mailto:paffenholz@opt.tu-darmstadt.de}{paffenholz@opt.tu-darmstadt.de}}

\address[Thomas Rehn]{
 initOS GmbH \& Co.\ KG,
 Hegelstr.\ 28,
 39104 Magdeburg, Germany
}
\email{\href{mailto:thomas.rehn@research.initos.com}{thomas.rehn@research.initos.com}}

\thanks{M.~Joswig, B.~Lorenz and A.~Paffenholz are partially supported by the DFG within the Priority Program 1489. M.~Joswig is additionally supported by Einstein Foundation Berlin.}


\begin{document}
\begin{abstract}
  The main purpose of this paper is to report on the state of the art of computing integer hulls and their facets as
  well as counting lattice points in convex polytopes.  Using the \polymake system we explore various algorithms and
  implementations.  Our experience in this area is summarized in ten ``rules of thumb''.
\end{abstract}
\maketitle

\section{Introduction}

\noindent
In integer and linear optimization the software workhorses are solvers for linear programs (based on simplex or interior
point methods) as well as generic frameworks for branch-and-bound or branch-and-cut schemes.  Comprehensive
implementations are available both as Open Source, like \scip \cite{scip}, as well as commercial software, like \cplex
\cite{cplex} and \gurobi \cite{gurobi}.  While today it is common to solve linear programs with millions of rows and
columns and, moreover, mixed integer linear programs with sometimes hundreds of thousands of rows and columns, big
challenges remain.  For instance, the 2010 version of the \miplib \cite{miplib} lists the mixed-integer problem
\texttt{liu} with $2\thinspace178$ rows, $1\thinspace156$ columns, and a total of only $10\thinspace626$ non-zero
coefficients; this seems to be impossible to solve with current techniques.  One way to make progress in the field is to
invent new families of cutting planes, either of a general kind or specifically tailored to a class of examples.  In the
latter situation the strongest possible cuts are obviously those arising from the facets of the (mixed) integer hull.  A main
purpose of this note is to report on the state of the art of getting at such facets in a brute force kind of way.  And
we will do so by explaining how our software system \polymake~\cite{polymake_prog} can help.

Here we focus on integer linear programming (ILP); mixed integer linear programming (MILP) will only be mentioned in
passing.  To avoid technical ramifications we assume that all our linear programs (LP) are bounded.  The brute force
method for obtaining all facets of the integer hull is plain and simple, and it has two steps.  First, we compute all the feasible integer points.  Since we assumed
boundedness these are only finitely many.  Second, we compute the facets of their convex hull.  Of course, the catch is
that neither problem is really easy.  Deciding if an ILP has an (integer) feasible point is known to be NP-complete
\cite{GareyJohnson79}.  So, we may not even hope for any efficient algorithm for the first step.  Most likely, the
situation for the second is about equally bad.  While it is open whether or not there is a convex hull algorithm which
runs in polynomial time measured in the combined sizes of the input and the output, recent work of Khachiyan \emph{et
al.}~\cite{KhachiyanEtAl08} indicates a negative answer.  They show that computing the vertices of an unbounded
polyhedron is hard; the difference to the general convex hull problem is that their result does not say anything about
the rays of the polyhedron.

Our paper is organized as follows.  We start out with a very brief introduction to the \polymake system and its usage.
In Section~\ref{sec:convex_hull} we explore how various convex hull algorithms and
their implementations behave on various kinds of input.  Our input is chosen according to typical scenarios motivated
by questions in optimization.  Deliberately we picked data within a manageable range; our goal is
to give a feel for what can be done within about one hour of CPU time on a current standard desktop machine.
For particularly hard convex hull problems see \cite{AvisBremnerSeidel97} and~\cite{Beneath-Beyond}.
Section~\ref{sec:integer_points} is devoted to
enumerating lattice points in polytopes, which is actually the first step in the procedure sketched above.  Like for the
convex hull computations there are several methods which behave quite differently depending on the input.  We would like
to stress that all the software systems mentioned have their preferred types of input for which they are (often vastly)
superior to the others.  There are no globally optimal algorithms known, neither for convex hull computations nor for
lattice point enumeration.  Even worse, in general, it is very difficult to say which method works best on which input.
Only a thorough  geometric and combinatorial analysis a posteriori allows for precise
statements.  Yet we will try to sum up our experience in this area in several ``rules of thumb'', all of which have to
be taken with a grain of salt.
The software systems tested are interfaced to by \polymake, and some of them are even shipped with it.
All our experiments are run through \polymake, which does create very little overhead.
The current version of \polymake does not directly make use of our rules of thumb.
Instead the user needs to pick the best method according to her or his judgment (or rely on some default).
We close the paper with additional references to the literature and related experimental results.
See the \hyperref[sec:setup]{Appendix} for details on our experimental setup.

Computers of today and the foreseeable future are faster than their predecessors just because they provide more cores.
Hence the parallelization of algorithms and software will play an increasingly important role.
Yet here we largely restrict our attention to single-threaded computations in order to not stretch our ---already quite involved--- test scenarios beyond any manageable limits.
As an exception we report on one parallel test with \normaliz \cite{normaliz2}.
The very recent version 6.0 of \lrs \cite{lrslib,lrs-paper} offers single machine shared memory multicore parallelization as well as \texttt{MPI} based multi-machine parallelization.
See Avis and Jordan \cite{AvisJordan:1510.02545} for a computational test of several large examples.

\subsubsection*{Acknowledgment}
We thank Thomas Opfer for contributing his implementation of the dual simplex method, originally written
for his Master's Thesis~\cite{ThomasDipl}, to the \polymake project, and this includes the code maintenance until today.

Moreover, we are very grateful to the developers of \cdd, \lrs, \normaliz and \ppl as they gave us various kind of valuable feedback.
The comments by David Avis and Winfried Bruns were particularly detailed.
What we found most rewarding is the fact that, partially in reaction to the 2014 preprint version of this paper, the teams of \lrs and \normaliz published new releases of their codes.
Throughout these show improvements which are sometimes very substantial, e.g., for \normaliz' handling of non-symmetric cut polytopes.
The interested reader may find it worth-while to compare the results below with that preprint version (which is still available as \href{http://arxiv.org/abs/1408.4653v1}{arXiv:1408.4653v1}).

\section{Example \polymake session}
\label{sec:session}

\noindent
In this section we will give a few examples on how one can use \polymake for various tasks, \textit{e.g.}, convex and
integer hull computations or linear programming.  The main application area of \polymake are polytopes and polyhedra,
but the software can also deal with a number of other kinds of objects. That is why \polymake is split into different
so-called \emph{applications}, for instance, for graphs, fans, ideals, matroids, simplicial complexes, or toric
varieties. We will only work with the application polytope here, which is indicated by the corresponding prompt
\lstinline|polytope>| in front of each code line.  The \polymake system employs a rule based mechanism for deciding
which algorithms satisfy user requests; this is explained in~\cite{DMV:polymake}. However, one may also force the usage
of specific algorithms. Thus, a big advantage of \polymake is that the user can use and compare various algorithms through
a common interface.  Our special object model, which here is used for the polytopes and the linear programs, is the
topic of~\cite{Flexible+Object+Hierachies+in+polymake}. The software is available from
\href{http://www.polymake.org}{polymake.org}~\cite{polymake_prog}.

\medskip

For our first example, a \emph{fractional knapsack polytope} $P$ is the set of points in the non-negative orthant of, say $\RR^d$, which
satisfies one extra linear inequality, \textit{i.e.},
\begin{align}
  P\ =\ \left\{\; x\in \RR^d\mid \left\langle a,x\right\rangle\,\le\, b\,,\; x_i\,\ge\, 0\quad\text{for all $i$}\right\}\label{eq:knapsack}
\end{align}
for some integral normal vector $a\in \ZZ^d$ and $b\in \ZZ$.  Such polytopes are easy to construct in \polymake.  For
instance, the command
\begin{lstlisting}
polytope> $k = fractional_knapsack([40,-2,-3,-5,-8,-13]);
\end{lstlisting}
%$
produces the fractional knapsack polytope for the inequality
\[
2 x_1 + 3 x_2 + 5 x_3 + 8 x_4 + 13 x_5 \ \le \ 40
\]
in $\RR^5$.  If all coefficients are positive, as in our example, the polytope is bounded.  In this case the combinatorics
is not so exciting, as all fractional knapsack polytopes are simplices.  Each line in the following output corresponds to one of
the six vertices.
\begin{lstlisting}
polytope> print $k->VERTICES;
1 0 0 0 0 0
1 20 0 0 0 0
1 0 40/3 0 0 0
1 0 0 8 0 0
1 0 0 0 5 0
1 0 0 0 0 40/13
\end{lstlisting}
%$
First, note that the coordinates contain rational numbers as \polymake uses \gmp's arbitrary precision rational
arithmetic \cite{gmp} by default.
Second, notice the leading $1$ of each vertex. This is due to the fact that polyhedra in \polymake are modeled as the intersection of
a cone with the affine hyperplane defined by $x_0=1$. With this homogenization many algorithms can deal with unbounded
polyhedra without much extra effort; e.g., see \cite[\S3.4]{Polyhedral+and+Algebraic+Methods}.

At this point \polymake already used a (dual) convex hull algorithm to find the vertices of the polytope which was given
by an outer description. There are several convex hull implementations to choose from within \polymake; see
Section~\ref{sec:convex_hull} for more details.  This particular fractional knapsack polytope belongs to the
\emph{Fibonacci fractional knapsack polytopes}.  Their integer hulls, the \emph{Fibonacci knapsack polytopes},
will be investigated in Section~\ref{sec:knapsack:hull} below.

\smallskip

Now we want to present how to optimize over this polytope with \polymake. Like for the convex hull algorithms, there are various (simplex type) implementations for linear programming available.
The first step is to create a \texttt{LinearProgram} object with an objective function. Then we can ask for the
maximal/minimal value and an optimal solution.
\begin{lstlisting}
polytope> $lp = $k->LP(LINEAR_OBJECTIVE=>new Vector([0,1,2,1,2,1]));
polytope> print $lp->MAXIMAL_VERTEX;
1 0 40/3 0 0 0
polytope> print $lp->MAXIMAL_VALUE;
80/3
\end{lstlisting}
Observe that the maximal vertex is not integral. To find an integral optimal solution you can repeat this process with
the integer hull of that polytope.

%polytope> print $k->N_LATTICE_POINTS;
%1366
\begin{lstlisting}
polytope> $ik = integer_hull($k);
polytope> print $ik->N_POINTS;
1366
polytope> print $ik->N_VERTICES;
16
\end{lstlisting}%$
Here \polymake enumerated all integral points inside the polytope. For this task there are again several different
algorithms available, see Section~\ref{sec:integer_points}. Afterwards we define a new polytope with all
those integral points as input points, and which is stored in \texttt{\$ik}.  Notice that \polymake also provides a
function \texttt{mixed\_integer\_hull} with the obvious interpretation.

\begin{figure}[t]
  \centering
  \subcaptionbox{Visualization of the integer points of the projection of \texttt{\$k} onto the first three coordinates with a maximal linear and a maximal integer
  solution.\label{fig:fibonacci_lattice_visual}}{\input{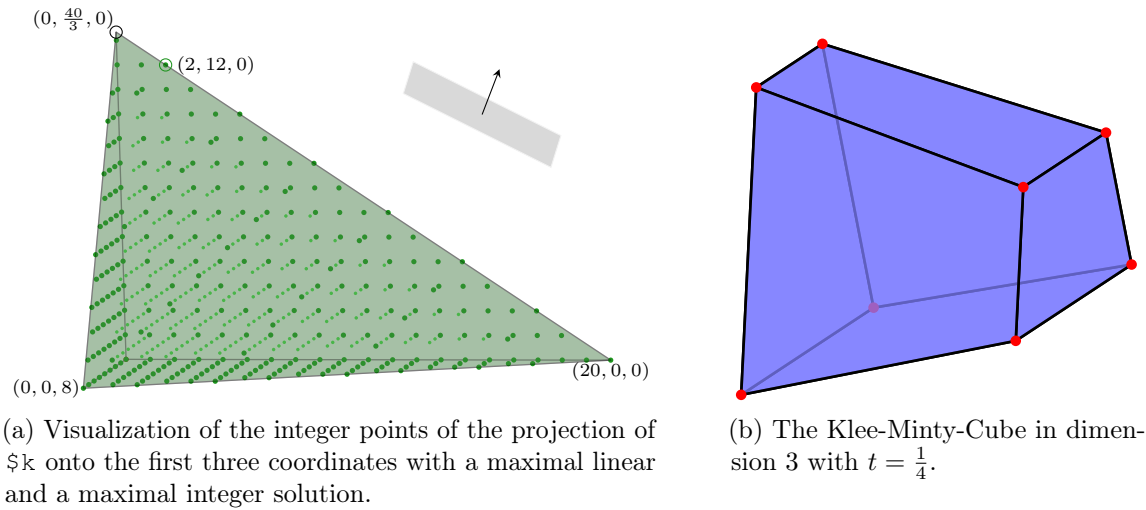}} \hfill
  \subcaptionbox{The
  Klee-Minty-Cube in dimension~$3$ with $t=\frac{1}{4}$.\label{fig:kmc}}{% polymake for assarf
% Mon Oct 26 10:52:53 2015
% unnamed
\begin{tikzpicture}[x  = {(0.9cm,-0.076cm)},
                    y  = {(-0.06cm,0.95cm)},
                    z  = {(-0.44cm,-0.29cm)},
                    scale = 2,
                    color = {lightgray}]

  % POINTS STYLE
  \definecolor{pointcolor_unnamed__1}{rgb}{ 1,0,0 }
  \tikzstyle{pointstyle_unnamed__1} = [fill=pointcolor_unnamed__1]

  % DEF POINTS
  \coordinate (v0_unnamed__1) at (2, 0.5, 1.875);
  \coordinate (v1_unnamed__1) at (2, 1.5, 1.625);
  \coordinate (v2_unnamed__1) at (2, 1.5, 0.375);
  \coordinate (v3_unnamed__1) at (2, 0.5, 0.125);
  \coordinate (v4_unnamed__1) at (0, 0, 0);
  \coordinate (v5_unnamed__1) at (0, 2, 0.5);
  \coordinate (v6_unnamed__1) at (0, 2, 1.5);
  \coordinate (v7_unnamed__1) at (0, 0, 2);

  % EDGES STYLE
  \definecolor{edgecolor_unnamed__1}{rgb}{ 0,0,0 }

  % FACES STYLE
  \definecolor{facetcolor_unnamed__1}{rgb}{ 0.5,0.5,1 }

  \tikzstyle{facestyle_unnamed__1} = [fill=facetcolor_unnamed__1, fill opacity=0.75, draw=edgecolor_unnamed__1, line width=1 pt, line cap=round, line join=round]

  % FACES and EDGES and POINTS in the right order
  \draw[facestyle_unnamed__1] (v0_unnamed__1) -- (v7_unnamed__1) -- (v4_unnamed__1) -- (v3_unnamed__1) -- (v0_unnamed__1) -- cycle;
  \draw[facestyle_unnamed__1] (v4_unnamed__1) -- (v7_unnamed__1) -- (v6_unnamed__1) -- (v5_unnamed__1) -- (v4_unnamed__1) -- cycle;
  \draw[facestyle_unnamed__1] (v4_unnamed__1) -- (v5_unnamed__1) -- (v2_unnamed__1) -- (v3_unnamed__1) -- (v4_unnamed__1) -- cycle;

  %POINTS
  \fill[pointcolor_unnamed__1] (v4_unnamed__1) circle (1 pt);

  %FACETS
  \draw[facestyle_unnamed__1] (v5_unnamed__1) -- (v6_unnamed__1) -- (v1_unnamed__1) -- (v2_unnamed__1) -- (v5_unnamed__1) -- cycle;

  %POINTS
  \fill[pointcolor_unnamed__1] (v5_unnamed__1) circle (1 pt);

  %FACETS
  \draw[facestyle_unnamed__1] (v1_unnamed__1) -- (v0_unnamed__1) -- (v3_unnamed__1) -- (v2_unnamed__1) -- (v1_unnamed__1) -- cycle;

  %POINTS
  \fill[pointcolor_unnamed__1] (v3_unnamed__1) circle (1 pt);
  \fill[pointcolor_unnamed__1] (v2_unnamed__1) circle (1 pt);

  %FACETS
  \draw[facestyle_unnamed__1] (v6_unnamed__1) -- (v7_unnamed__1) -- (v0_unnamed__1) -- (v1_unnamed__1) -- (v6_unnamed__1) -- cycle;

  %POINTS
  \fill[pointcolor_unnamed__1] (v6_unnamed__1) circle (1 pt);
  \fill[pointcolor_unnamed__1] (v7_unnamed__1) circle (1 pt);
  \fill[pointcolor_unnamed__1] (v0_unnamed__1) circle (1 pt);
  \fill[pointcolor_unnamed__1] (v1_unnamed__1) circle (1 pt);

  %FACETS

\end{tikzpicture}}
  \caption{Visualization of the two example polytopes, made with \polymake's interface to TikZ~\cite{TikZ}.}
\end{figure}

Solving an integer linear program is the same as doing linear optimization over the integer hull.  So we can proceed as before.
A projection of this Fibonacci polytope, the integral points, and the optimal vertices is shown in
Figure~\ref{fig:fibonacci_lattice_visual}.
\begin{lstlisting}
polytope> $ilp = $ik->LP(LINEAR_OBJECTIVE=>new Vector([0,1,2,1,2,1]));
polytope> print $ilp->MAXIMAL_VERTEX;
1 2 12 0 0 0
polytope> print $ilp->MAXIMAL_VALUE;
26
\end{lstlisting}

As already mentioned, \polymake usually works with exact rational coordinates.
However, it is also possible to work with other coordinate types, which includes certain extensions of rationals.
For instance, fields of Puiseux fractions are relevant in optimization since they can be used to model coordinate perturbations.
The \emph{(rational) Puiseux fractions} form a field whose elements are univariate rational functions with rational coefficients and rational exponents; they can be equipped with a natural ordering.
Computing with Puiseux fractions amounts to computing with the usual rational numbers and an infinitesimally small transcendental number; see \cite{MACIS2015} for the details.
We illustrate this with the $3$-dimensional \emph{Klee--Minty} cube~\cite{KleeMinty:1972} which can be constructed as follows.
\begin{lstlisting}
polytope> $monomial=new UniMonomial<Rational,Rational>(1);
polytope> $t=new PuiseuxFraction<Min>($monomial);
polytope> $c = klee_minty_cube(3,$t);
polytope> print_constraints($c);
Facets:
0: x1 >= (0)
1: -x1 >= (- 1)
2: -(t) x1 + x2 >= (0)
3: -(t) x1 - x2 >= (- 1)
4: -(t) x2 + x3 >= (0)
5: -(t) x2 - x3 >= (- 1)
\end{lstlisting}
The parentheses indicate that the coefficients are Puiseux fractions, rather than rational numbers.
Yet it makes sense to contemplate evaluating these expressions for special values.
For instance, $t=0$ yields the ordinary $0/1$-cube; and for $0<t<1/2$ we obtain the Klee--Minty cube.
Notice that, by construction, the latter inequality always holds in the field of Puiseux fractions.
As a benefit from working over Puiseux fractions, e.g., we can obtain the volume as a rational function in $t$ (which is a polynomial in this case):
\begin{lstlisting}
polytope> print $c->VOLUME;
(1 -2*t + t^2)
\end{lstlisting}%$

\medskip

In addition to the various applications shipped with \polymake, there are also third-party extensions, which add
further functionality. Of particular interest in the context of integer linear programming is the
\texttt{unimodularity-test} extension by Matthias Walter~\cite{unimodularity-test}, which implements an algorithm of Walter and Tr\"umper \cite{WalterTruemper:2013}.  Our wiki
page at \url{polymake.org} contains several tutorials, some of which specialize on topics in optimization.  In
particular, it features a proof-of-concept implementation of the branch-and-bound algorithm, to solve integer linear
programs without relying on computing all lattice points within the feasible domain, which is useful for educational
purposes.  Via the interface to \normaliz~\cite{normaliz2} our software can also compute Hilbert bases.  For instance,
this allows for computing Gomory-Chv\'atal closures by calling the function \texttt{gc\_closure}.

\medskip

\polymake also allows to directly read and write data files in the classical \texttt{lp}-format for linear programs and in \porta
format for convex hull problems in primal and dual form. This can be done with the commands \texttt{lp2poly, poly2lp,
  poly2porta}, and \texttt{porta2poly}. See the corresponding help text for the correct calling syntax and options
(e.g., \texttt{help "lp2poly";} inside the \polymake shell).

\medskip

The programming language of the \polymake shell is essentially \perl.  However, we extended the language in several
ways.  The most important modification is related to our implementation of an interface between \perl and \cplus, which
is based on \perl's XS interface.  For instance, calling the function \verb+new Vector(+$\ldots$\verb+)+ creates an
object via \cplus code, and the \perl shell keeps a reference to an opaque object.  The \cplus code has a template
parameter for the coefficients which defaults to \texttt{Rational}, which is our wrapper to the \gmp.  We modified \perl
to allow for template parameters in the \cplus data type.  For instance, \verb+new Vector<Integer>(+$\ldots$\verb+)+ is
legal in \polymake.  New instantiations of data types due to an uncommon choice of template parameters trigger the
automatic compilation of the appropriate \cplus code.  That code is kept such that it does not need to be recompiled
again.

\section{Convex Hull Computations}
\label{sec:convex_hull}
% d: dimension
% n: number of vertices
% m: number of facets

\noindent
The \emph{convex hull problem} in its primal form asks to compute the facets of the convex hull of finitely many given
points in $\RR^d$.  By cone polarity this is equivalent to computing the vertices and rays of a polyhedron which is
given in terms of finitely many linear equations and inequalities.  The latter form is the \emph{dual convex hull
problem}.  Both forms occur naturally in the context of linear and integer programming: \textit{e.g.},
to compute the
facets of the integer hull of a given polytope is a primal convex hull problem.  Dual convex hull problems arise,
\textit{e.g.}, in the computation of Voronoi diagrams.  Out of the four types of convex hull problems we will investigate two, namely computing the facets of integer hulls and computing Voronoi
diagrams.
Other applications of convex hull algorithms aim at certain parameterized optimization problems.
An example is the search for biologically correct alignments of chicken hemoglobin mRNA data \cite{ASCB:DeweyWoods}.

\subsection{Complexity Status}

The precise complexity status of the convex hull problem is unsettled: It is not known whether or not there exists an
algorithm whose complexity is bounded by a polynomial in the combined size of the input and the output. This concept is
called \emph{polynomial total time}; see~\cite{Valiant79}.  As a measure for the complexity of an algorithm it differs
from what is common in theoretical computer science, where algorithms are typically measured in the size of the input
only.  To evaluate convex hull algorithms from a practical point of view restricting the attention to the input size
only is somewhat unsatisfying as the size of the output can vary over several orders of magnitude.  To give just one
basic example: a cube in dimension $d$ has $2^d$ vertices but only $2d$ facets, while the dual cross polytope has $2d$
vertices but as many as $2^d$ facets.  By McMullen's Upper Bound Theorem \cite[\S8]{Ziegler95} the number of facets of a
$d$-polytope with $n$ vertices is bounded by $O(n^{d/2})$, and this bound is tight.  For fixed dimension $d$ the bound
$O(n^{d/2})$ becomes polynomial in $n$, and for this setup Chazelle gave a fixed dimension polynomial time algorithm
which is worst-case optimal with respect to the size of the input \cite{Chazelle93}.
However, by taking the dimension as a constant, such an algorithm is allowed to ``waste'' an exponential amount of time on many small convex hull problems.

More recently, Khachiyan \textit{et al.}
showed that it is $\#$P-hard to enumerate the \emph{vertices} of an unbounded
polyhedron given in terms of inequalities \cite{KhachiyanEtAl08}.  However, this result does not say anything about the
\emph{rays} of that polyhedron, whence it does not have a direct implication for the complexity status of the convex
hull problem.  Yet this is a strong hint that a polynomial total time convex hull algorithm does not exist.

\subsection{Common Algorithms}

In practice several algorithms and their implementations are used; see \cite{AvisBremnerSeidel97} and
\cite{Beneath-Beyond} for previous computational studies.  One purpose of this paper is to point out that essentially
all of them are relevant for various convex hull problems arising in linear and integer optimization.

The \emph{beneath-and-beyond} method is best described as a primal convex hull algorithm; it occurs in work of
Seidel~\cite{Seidel81}, see also \cite{Edelsbrunner87} and~\cite{Beneath-Beyond}.  So the goal is to compute the facets
of the convex hull, say $P$, of finitely many given points.  We may assume that the points affinely span the entire
space for otherwise we can restrict our attention to the affine span.  Picking an affine basis among the input points
yields a simplex $P_{1}$ of full dimension $d$
which is contained in $P$.  By solving systems of linear equations we
obtain the facets of $P_{1}$.  Now the remaining input points are added inductively, yielding a sequence of polytopes
$P_1,P_2,\dots,P_{n-d}=P$.  In the inductive step triangulations of the polytopes $P_i$ are maintained, and these allow
for the computations of the facets.  The size of the final triangulation of the polytope $P$ is the decisive factor in
the running time; see \cite[Proposition 3.2]{Beneath-Beyond}.  Notice that there are polytopes with few vertices and few
facets for which nonetheless all triangulations are large; see \cite{AvisBremnerSeidel97} and
\cite[\S3.2]{Beneath-Beyond}.  In our tests we use the implementation included in \polymake,
which is marked as ``\bb'' below.

Algorithms which inductively compute the intermediate convex hulls $P_i$ for all $i$ are called
\emph{incremental}. Bremner showed that no incremental convex hull algorithm can have a polynomial total running time
\cite{Bremner99}.

A dual method, which is also incremental, is \emph{double description}. It was first described by Motzkin \emph{et
  al.}~\cite{Motzkin53}; see also \cite{FukudaProdon96} and \cite[\S5.2]{Polyhedral+and+Algebraic+Methods}.
This is essentially the same as \emph{Fourier--Motzkin
  elimination} and \emph{Chernikova's algorithm} \cite{Chernikova64}.  Given finitely many inequalities, one starts with
a subset of the inequalities which enclose a full-dimensional simplex, $Q_1$.  Up to solving a linear programming
feasibility problem and up to a projective transformation, this is not a restriction.  The vertices of that simplex are
again computed by solving systems of linear equations.  In the inductive step, to obtain $Q_{i+1}$, the hyperplane
defined by a new inequality separates the vertices of $Q_i$ into those which are also vertices of $Q_{i+1}$ and those
which are not.  The new vertices of $Q_{i+1}$ arise as intersection of that hyperplane with the edges of $Q_i$.  The
extra data structure which needs to be maintained is the vertex-edge graph of the intermediate polytopes.  The double
description method is implemented, \textit{e.g.}, in \cdd \cite{cddlib,cdd-paper}, \ppl \cite{ppl} and \porta \cite{porta}.  Those software
packages are tested below. See \cite{Borgwardt07} for an average-case analysis of the beneath-and-beyond and double description methods.

From a certain point of view one can even consider beneath-and-beyond and double description as similar;
\emph{e.g.}, see \cite{FukudaProdon96}.  However, there is a choice which kind of data structures to maintain and how
to maintain them during the computation.  In our terminology, beneath-and-beyond primarily uses a triangulation (in the primal
setting), while double description primarily uses the graph (in the dual setting).  As our experiments below show
this leads to a very different behavior on various kinds of input.

A third convex hull algorithm which we consider here is \emph{reverse search} by Avis and Fukuda~\cite{AvisFukuda96}.
The only non-incremental algorithm in our selection is based on the simplex method for linear programming.  Any linear
objective function which is in general position with respect to a given polytope $R$, induces a direction on each edge
of $R$.  Starting from some vertex the optimum (which is unique as the objective function is generic) can be reached by
following a directed path.
The reverse search algorithm in its natural form computes dual convex hulls. It solves
one linear program to determine one vertex of the feasible region from the given list of finitely many linear
inequalities.  By pivoting backwards all vertices can be traced.  The non-trivial ingredient in this algorithm is that
this can be accomplished with constant memory after the initialization.  The reverse search method is implemented in
\lrs \cite{lrslib,lrs-paper}. There is a multithreaded version of \lrs, which we did not test; see \cite{AvisJordan:1510.02545}.

There are more convex hull algorithms and implementations. Bruns, Ichim and S\"oger suggested an interesting hybridization of the beneath-and-beyond method with Fourier--Motzkin elimination, which they called \emph{pyramid decomposition} \cite{BrunsIchimSoeger:1206.1916}.
As a special feature their \normaliz implementation supports multithreading; we report on some tests with \normaliz, and these are marked `\texttt{nmz}'.
The \cgal library \cite{cgal} also offers convex hull implementations, but we did not test them here; the reason is that the output is somewhat different.
Joswig and Ziegler described a convex hull algorithm based on simplicial homology computation \cite{JoswigZiegler04}, which is interesting from a theoretical point of view only.

Our first two rules of thumb serve as a first approximation, and this is consistent with what had previously been observed by others.
\begin{ruleofthumb}\label{rule:double-description}
  If you do not know anything about your input to a convex hull problem, try double description.
\end{ruleofthumb}
While this is useful to keep in mind, a large part of this paper is devoted to discuss \emph{natural} situations in which this simple minded strategy leads astray.
A conceptual disadvantage of double description, being an incremental algorithm, is that its entire output is only available at the very end of the computation.
Non-incremental techniques, such as reverse search, do not share this disadvantage.
This immediately leads to the next rule of thumb.
\begin{ruleofthumb}
  Use reverse search if you expect the output to be extremely large, and if partial information is useful.
\end{ruleofthumb}
\lrs features options for obtaining an unbiased estimate of the output size (and \lrs running time), from looking into the first few levels of the search tree.
By comparing the number of vertices (or facets) with the number of bases one can also get an idea of how degenerate the input is.
Notice, however, that \polymake's current interface does not support these advanced features of \lrs.

\subsection{Non-Symmetric Cut Polytopes}
\label{sec:non-symmetric-cut}

Let us start out with a class of examples where the first rule above works to our advantage.  Consider a finite simple
graph $G=(V,E)$ which is undirected, \textit{i.e.}, the \emph{nodes} form a finite set $V$
and the set $E$ of edges consists
of two-element subsets of $V$.  A (possibly trivial) partition $V=A+B$ of the node set induces a \emph{cut}
$C(A,B):=\{e\in E\mid |e\cap A|=|e\cap B|=1\}$ in $G$. We can encode such a cut via its \emph{incidence vector}
$\chi^C\in\{0,1\}^{|E|}$. The \emph{cut polytope} of $G$ is the polytope
\[
\Cut(G)\ :=\ \conv(\chi^C\mid C\text{ is a cut in G})\,.
\]
A classical line of research in combinatorial optimization asks to determine the facets of cut polytopes; see
Schrijver~\cite[\S75.7]{Schrijver03:CO_C}.

\begin{remark}\label{rem:symmetric-cut}
  Each automorphism of $G$ gives rise to a linear automorphism of its cut polytope $\Cut(G)$.
\end{remark}

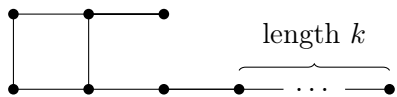
\begin{figure}[htb]
  \centering
  \begin{tikzpicture}[scale=1]
  \tikzstyle{edge} = [draw,black,-]

  \coordinate (0) at (1,1);
  \coordinate (1) at (0,1);
  \coordinate (2) at (-1,1);
  \coordinate (3) at (-1,0);
  \coordinate (4) at (0,0);
  \coordinate (5) at (1,0);
  \coordinate (6) at (2,0);
  \coordinate (7) at (4,0);
  \coordinate (dots) at (3,0);
  
  \coordinate (bracestart) at (2,0.25);
  \coordinate (braceend) at (4,0.25);

  \draw[edge] (0) -- (1) -- (2) -- (3) -- (4) -- (1) -- cycle;

  \draw[edge] (4) -- (5);
  \draw[edge] (5) -- (6);
  \draw[edge] (5) -- (dots);
  \draw[edge] (dots) -- (7);
  
  \node[circle, fill=white] at (dots) {$\dots$};
  
  \foreach \point in {0,1,2,3,4,5,6,7}
    \fill[black] (\point) circle (2pt);
    
  \draw [decoration={brace}, decorate] (bracestart) -- (braceend);
  \node[above] at (3,0.4) {length $k$};
\end{tikzpicture} 
  \caption{Asymmetric graph $G_k$ on $k+6$ nodes.}
  \label{fig:graph}
\end{figure}
For our first experiment we look at graphs with few nodes and which are completely asymmetric, \textit{i.e.}, their
automorphism group is trivial.  In Section~\ref{sec:symmetric-cut} below we will discuss cut polytopes of graphs which
\emph{are} symmetric.  Let $G_k$ be the connected graph with $k+6$ nodes and $k+6$ edges shown in
Figure~\ref{fig:graph}.  For $k\ge 1$ the graph $G_k$ is completely asymmetric, while $G_0$ admits one involutory
automorphism.

The relevant parameters for the cut polytopes of the graph $G_k$ are easy to determine.  The dimension,
$d=\dim(\Cut(G_k))$, equals the number of edges, which is $k+6$.  There are $n=2^{k+5}$ many cuts (including the empty
cut), all of which form vertices, as cut polytopes are $0/1$-polytopes.
The values $d$ and $n$ describe the size of the
input.  In this case, since the graph $G_k$ does not contain a $K_5$-minor, the facets of $\Cut(G_k)$ can be determined
via a result of Barahona~\cite{Barahona83}; see also \cite[Corollary~75.4f]{Schrijver03:CO_C}.  This yields
$m=2d+8=2k+20$ facets.

Our experiment investigates the primal convex hull computations for the polytopes $\Cut(G_k)$, where $k$ ranges from $0$
to $14$, and for which we tried the six implementations \bb, \cdd, \lrs, \nmz, \porta, and \ppl. The polytopes were
given to \polymake in reduced vertex description via \texttt{VERTICES} and \texttt{LINEALITY\_SPACE}. See the accompanying
web page~\cite{polymake_exp} for full details.

\polymake has default options for computing convex hulls etc.\ but you can directly influence the choice of
algorithm for various tasks within the \polymake shell by issuing \texttt{prefer} statements. To choose, e.g., \ppl as
the default method for convex hull computations you can use
\begin{lstlisting}
polytope> prefer "ppl";
\end{lstlisting}
This choice will be stored in the preferences file of \polymake (\texttt{.polymake/prefer.pl} in your home folder) and is
valid until you set another default. To switch only temporarily you can use \texttt{prefer\_now}, which is only valid on
the line it is contained in, e.g.,
\begin{lstlisting}
polytope> prefer_now "lrs"; print $p->FACETS;
\end{lstlisting}%$
This can be fine tuned with minor keys in the form \texttt{ppl.convex\_hull}, see the \polymake wiki for more details.
By the way, omitting the ``\texttt{print}'' in the second command above will trigger the convex hull computation without producing output.
This is instrumental, e.g., for proper timings.

\begin{figure}[tb]
  \centering
  % This file was created by matlab2tikz v0.4.5 running on MATLAB 8.2.
% Copyright (c) 2008--2014, Nico Schlömer <nico.schloemer@gmail.com>
% All rights reserved.
% 
% The latest updates can be retrieved from
%   http://www.mathworks.com/matlabcentral/fileexchange/22022-matlab2tikz
% where you can also make suggestions and rate matlab2tikz.
% 
%
% defining custom colors
\definecolor{mycolor1}{rgb}{0.16863,0.50588,0.33725}%
\begin{tikzpicture}

\begin{axis}[%
width=6.4cm,
height=5cm,
scale only axis,
xmin=0,
xmax=14,
ymode=log,
ymin=0.001,
ymax=12000,
yminorticks=true,
legend style={at={(0.71,0.02)},anchor=south west,draw=black,fill=white,legend cell align=left}
]
\addlegendimage{bbstyle};
\addlegendentry{\bb};
\addlegendimage{cddstyle};
\addlegendentry{\cdd};
\addlegendimage{lrsstyle};
\addlegendentry{\lrs};
\addlegendimage{nmzstyle};
\addlegendentry{\nmz\!\!$_6$};
\addlegendimage{portastyle};
\addlegendentry{\porta};
\addlegendimage{pplstyle};
\addlegendentry{\ppl};

\addplot [bbstyle]
  table[row sep=crcr]{
0	 0.008\\
1	 0.030\\
2	 0.168\\
3	 1.250\\
4 	14.411\\
};	
%\addlegendentry{\bb};

\addplot [cddstyle]
  table[row sep=crcr]{
0	  0.005\\
1	  0.015\\
2	  0.050\\
3	  0.184\\
4	  0.700\\
5	  2.877\\
6	 12.086\\
7	 50.064\\
8	210.074\\
9	974.681\\
};
%\addlegendentry{\cdd};

\addplot [lrsstyle]
  table[row sep=crcr]{
0	   0.005\\
1	   0.060\\
2	   0.937\\
3	  18.195\\
4	 329.093\\
5	7699.132\\
};
%\addlegendentry{\lrs};

\addplot [pplstyle]
  table[row sep=crcr]{
0       0.001\\
1       0.002\\
2       0.004\\
3       0.008\\
4       0.018\\
5       0.043\\
6       0.110\\
7       0.298\\
8       0.890\\
9       3.182\\
10     12.680\\
11     53.547\\
12    216.302\\
13   1148.847\\
14   3641.770\\
};
%\addlegendentry{\ppl};

\addplot [nmzstyle]
  table[row sep=crcr]{
0        0.001\\
1        0.002\\
2        0.005\\
3        0.010\\
4        0.021\\
5        0.047\\
6        0.114\\
7        0.293\\
8        0.841\\
9        2.799\\
10      10.233\\
11      42.461\\
12     181.790\\
13     952.263\\
14    4843.120\\
};
%\addlegendentry{\nmz\!\!$_6$};

\addplot [portastyle]
  table[row sep=crcr]{
0      0.065 \\
1      0.089 \\
2      0.107 \\
3      0.246 \\
4      0.916 \\
5      7.980 \\
6     98.987 \\
7   1810.960 \\
};
%\addlegendentry{\porta};

\end{axis}
\end{tikzpicture}%
  \caption{Timings (in seconds, logarithmically scaled)) for convex hull computations of non-symmetric cut polytopes $\Cut(G_k)$; see Table~\ref{tab:non_sym_cut}.  The curves for \texttt{nmz}$_6$ and \ppl almost agree.}
  \label{fig:non_sym_cut:timings}
\end{figure}
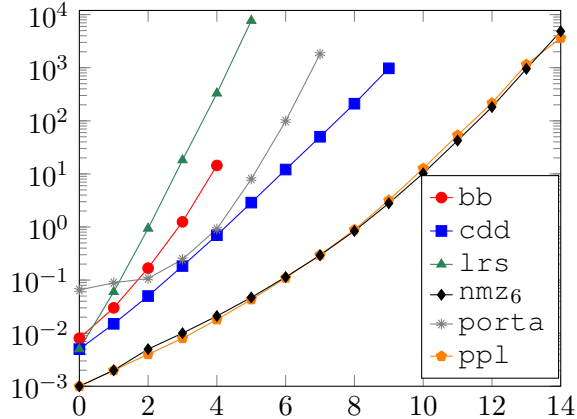

The programs \cdd, \lrs, \nmz and \ppl are interfaced via their library versions, hence there is virtually no overhead
from testing through the \polymake system.
However, our interface to \porta is based on an exchange of text files.
This does create some overhead, which is visible in \porta's curve in Figure~\ref{fig:non_sym_cut:timings}, for small input.
For large input, however, this effect can be neglected (at least on a logarithmic scale).

Every code was tried ten times on each input with the same parameters. The average running times are displayed in
Table~\ref{tab:non_sym_cut} and Figure~\ref{fig:non_sym_cut:timings}. Note that \bb runs out of memory quite early.
More precisely, our imposed memory limitation to 4GB does not suffice for $\Cut(G_5)$ and beyond; see the
\hyperref[sec:setup]{Appendix} for details on our experimental setup.  The reason is that the triangulations grow very
quickly.  For instance, the triangulation of $\Cut(G_4)$ which was computed by \bb already has
1\thinspace040\thinspace517 facets.  For $\Cut(G_5)$ we even get 11\thinspace199\thinspace900 facets within like 200s
but requiring about 10GB of main memory.
The asymptotic behavior of \cdd is quite good, but it is clearly outperformed by \ppl and \nmz.
The only parallel code is \normaliz, which here we tested with six cores.
Our experiments (undocumented here) suggest a factor of three in performance versus the single-threaded version of \nmz.

\begin{ruleofthumb}
  Do use double description for computing the facets of $0/1$-polytopes.
\end{ruleofthumb}

It is inspiring to see how well the curves for the six-threaded \nmz and the single-threaded \ppl match.
For practical purposes this seems to suggest that \nmz parallelization, which employs the \texttt{OpenMP} protocol, is particularly interesting for more than six cores.

\subsection{Knapsack Integer Hulls}
\label{sec:knapsack:hull}

For our second test we look at the kind of knapsack polytopes introduced in Section~\ref{sec:session}.  Consider the
Fibonacci sequence
\begin{align*}
  a_1 &= 2 & a_2&=3 & a_i &= a_{i-2} + a_{i-1} \,.
\end{align*}
Then the \emph{Fibonacci fractional knapsack polytope} in dimension $d$ for the parameter $b>0$ is defined as
\[
F_d(b) \ =\ \{ x \in \RR_{\ge0}^d \mid a^Tx \le b \} \,.
\]
As pointed out before, as mere polytopes these are boring as all of them are $d$-dimensional simplices.  What makes a
more interesting computational experiment is to look at the set of integer points in $F_d(b)$ and to compute their
convex hull, which yields the corresponding \emph{Fibonacci knapsack polytope}.  So, as in
Section~\ref{sec:non-symmetric-cut}, we are computing the facets of sets of integer points.  As a key difference,
however, here the vast majority of the input points is redundant, \textit{i.e.}, they do not form vertices of the
integer hull.  See Section~\ref{sec:knapsack:integer} below for a discussion about how to compute these integer
points. Each polytope for this set of experiments was defined in \emph{redundant} vertex description via the property
\texttt{POINTS}, see~\cite{polymake_exp} for details.

\begin{table}[t]
  \centering
  \caption{Properties of the Fibonacci fractional knapsack polytopes $F_d(b)$.}
  \label{tab:knapsack:properties}

  \begin{tabular}[t]{r@{\hspace{1cm}}rrrrrrr}
    \multicolumn{8}{c}{Number of integer points}\\
    \toprule
   $d \setminus b$   & \multicolumn{1}{c}{40} & \multicolumn{1}{c}{50} & \multicolumn{1}{c}{60} & \multicolumn{1}{c}{70} & \multicolumn{1}{c}{80} & \multicolumn{1}{c}{90} & \multicolumn{1}{c}{100}\\
    \midrule
      4 &  1021 &  2145 &  4008 &   6879 &  11069 &  16929 &  24853\\
      5 &  1366 &  3173 &  6509 &  12182 &  21245 &  35025 &  55157\\
      6 &  1481 &  3626 &  7853 &  15516 &  28544 &  49570 &  82090\\
    \bottomrule\\

    \multicolumn{8}{c}{Number of facets of the integer hull}\\
    \toprule
   $d \setminus b$   & \multicolumn{1}{c}{40} & \multicolumn{1}{c}{50} & \multicolumn{1}{c}{60} & \multicolumn{1}{c}{70} & \multicolumn{1}{c}{80} & \multicolumn{1}{c}{90} & \multicolumn{1}{c}{100}\\
    \midrule
    4 &   6 &   7 &   7 &    8 &   6 &    6 &    8 \\
    5 &  12 &  15 &  12 &   12 &   8 &   13 &   15 \\
    6 &  25 &  20 &  21 &   25 &  21 &   18 &   22 \\
    \bottomrule\\

    \multicolumn{8}{c}{Number of vertices of the integer hull}\\
    \toprule
   $d \setminus b$   & \multicolumn{1}{c}{40} & \multicolumn{1}{c}{50} & \multicolumn{1}{c}{60} & \multicolumn{1}{c}{70} & \multicolumn{1}{c}{80} & \multicolumn{1}{c}{90} & \multicolumn{1}{c}{100}\\
    \midrule
    4 &   8 &   11 &    9 &   12 &    8 &    8 &   12 \\
    5 &  16 &   25 &   19 &   23 &   13 &   19 &   25 \\
    6 &  35 &   37 &   35 &   40 &   35 &   31 &   40 \\
    \bottomrule

  \end{tabular}
  \bigskip
\end{table}

To get an idea about the input and output sizes Table~\ref{tab:knapsack:properties} lists the number of lattice points
(which provide the input), the number of facets (which form the output) and the number of vertices depending on $d$ and
$b$.  By a result of Hayes and Larman \cite{HayesLarman83} the number of vertices of the integer hull of any knapsack
$d$-polytope is bounded by $\sigma^d$, where $\sigma$ is the length of the binary encoding of the knapsack inequality;
our examples are rather far from this bound.

\begin{figure}[t]
  \centering
  \subcaptionbox{default order\label{fig:knapsack:timings:plot:normal}}{% This file was created by matlab2tikz v0.4.5 running on MATLAB 8.2.
% Copyright (c) 2008--2014, Nico Schlömer <nico.schloemer@gmail.com>
% All rights reserved.
% 
% The latest updates can be retrieved from
%   http://www.mathworks.com/matlabcentral/fileexchange/22022-matlab2tikz
% where you can also make suggestions and rate matlab2tikz.
% 
%
% defining custom colors
\definecolor{mycolor1}{rgb}{0.16863,0.50588,0.33725}%
\begin{tikzpicture}

\setlength{\mywidth}{.39\linewidth} % you might want to change this
\setlength{\myheight}{.8\mywidth} % do not change this

\begin{axis}[%
width=\mywidth, %6.4cm,
height=\myheight, %5cm,
scale only axis,
xmin=40,
xmax=100,
ymode=log,
ymin=0.01,
ymax=10000,
yminorticks=true,
legend style={at={(0.6,0.02)},anchor=south west,draw=black,fill=white,legend cell align=left, legend columns=2}
]

%\addlegendimage{bbstyle};
%\addlegendentry{\bb};
%\addlegendimage{cddstyle};
%\addlegendentry{\cdd};
%\addlegendimage{lrsstyle};
%\addlegendentry{\lrs};
%\addlegendimage{nmzstyle};
%\addlegendentry{\nmz};
%\addlegendimage{portastyle};
%\addlegendentry{\porta};
%\addlegendimage{pplstyle};
%\addlegendentry{\ppl};

\addplot [bbstyle]
  table[row sep=crcr]{
40        0.358\\
50        1.268\\
60        4.273\\
70       15.039\\
80       42.077\\
90      136.823\\
100     355.191\\
};	
%\addlegendentry{\bb};

\addplot [cddstyle]
  table[row sep=crcr]{
40         5.667 \\
50        30.601 \\
60       126.978 \\
70       447.216 \\
80      1352.025 \\
90      3594.319 \\
100     8737.240 \\
};
%\addlegendentry{\cdd};

\addplot [lrsstyle]
  table[row sep=crcr]{
40         27.849 \\
50        213.319 \\
60       1163.065 \\
70       4957.700 \\
};
%\addlegendentry{\lrs};

\addplot [nmzstyle]
  table[row sep=crcr]{
40       0.565\\
50       2.377\\
60       8.096\\
70      26.586\\
80      74.743\\
90     194.972\\
100    421.936\\
};
%\addlegendentry{\nmz};

\addplot [pplstyle]
  table[row sep=crcr]{
40      0.048 \\
50      0.123 \\
60      0.274 \\
70      0.657 \\
80      1.361 \\
90      3.472 \\
100     8.084 \\
};
%\addlegendentry{\ppl};

\addplot [portastyle]
  table[row sep=crcr]{
40      0.195  \\
50      0.685  \\
60      1.868  \\
70      5.363  \\
80     12.614  \\
90     37.579  \\
100    96.675  \\
};
%\addlegendentry{\porta};

\end{axis}
\end{tikzpicture}%
}
  \subcaptionbox{individual best order\label{fig:knapsack:timings:plot:best}}{% This file was created by matlab2tikz v0.4.5 running on MATLAB 8.2.
% Copyright (c) 2008--2014, Nico Schlömer <nico.schloemer@gmail.com>
% All rights reserved.
% 
% The latest updates can be retrieved from
%   http://www.mathworks.com/matlabcentral/fileexchange/22022-matlab2tikz
% where you can also make suggestions and rate matlab2tikz.
% 
%
% defining custom colors
\definecolor{mycolor1}{rgb}{0.16863,0.50588,0.33725}%
\begin{tikzpicture}

\setlength{\mywidth}{.39\linewidth} % you might want to change this
\setlength{\myheight}{.8\mywidth} % do not change this

\begin{axis}[%
width=\mywidth, %6.4cm,
height=\myheight, %5cm,
scale only axis,
xmin=40,
xmax=100,
ymode=log,
ymin=0.01,
ymax=10000,
yminorticks=true,
legend style={at={(0.4,0.02)},anchor=south west,draw=black,fill=white,legend cell align=left, legend columns=3}
]

\addlegendimage{bbstyle};
\addlegendentry{\bb};
\addlegendimage{cddstyle};
\addlegendentry{\cdd};
\addlegendimage{lrsstyle};
\addlegendentry{\lrs};
\addlegendimage{nmzstyle};
\addlegendentry{\nmz};
\addlegendimage{portastyle};
\addlegendentry{\porta};
\addlegendimage{pplstyle};
\addlegendentry{\ppl};

\addplot [bbstyle]
  table[row sep=crcr]{
40       0.054\\
50       0.140\\
60       0.363\\
70       1.188\\
80       3.457\\
90       9.689\\
100     23.175\\
};	
%\addlegendentry{\bb};

\addplot [cddstyle]
  table[row sep=crcr]{
40         5.667 \\
50        30.601 \\
60       126.978 \\
70       447.216 \\
80      1352.025 \\
90      3594.319 \\
100     8737.240 \\
};
%\addlegendentry{\cdd};

\addplot [lrsstyle]
  table[row sep=crcr]{
40         2.108 \\
50         9.002 \\
60        32.371 \\
70       111.586 \\
80       287.381 \\
90       599.524 \\
100      987.394 \\
};
%\addlegendentry{\lrs};

\addplot [nmzstyle]
  table[row sep=crcr]{
40       0.565\\
50       2.377\\
60       8.096\\
70      26.586\\
80      74.743\\
90     194.972\\
100    421.936\\
};
%\addlegendentry{\nmz};

\addplot [pplstyle]
  table[row sep=crcr]{
40      0.048 \\
50      0.123 \\
60      0.274 \\
70      0.657 \\
80      1.361 \\
90      3.472 \\
100     8.084 \\
};
%\addlegendentry{\ppl};

\addplot [portastyle]
  table[row sep=crcr]{
40      0.076 \\
50      0.213 \\
60      0.349 \\
70      0.685 \\
80      0.854 \\
90      2.278 \\
100     4.443 \\
};
%\addlegendentry{\porta};

\end{axis}
\end{tikzpicture}%
}
  \caption{Timings (in seconds) for integer hull computations of Fibonacci fractional knapsack polytopes $F_5(b)$ depending on $b$, on a logarithmic scale; see Table~\ref{tab:knapsack:timings} for exact timings.}
  \label{fig:knapsack:timings}
\end{figure}
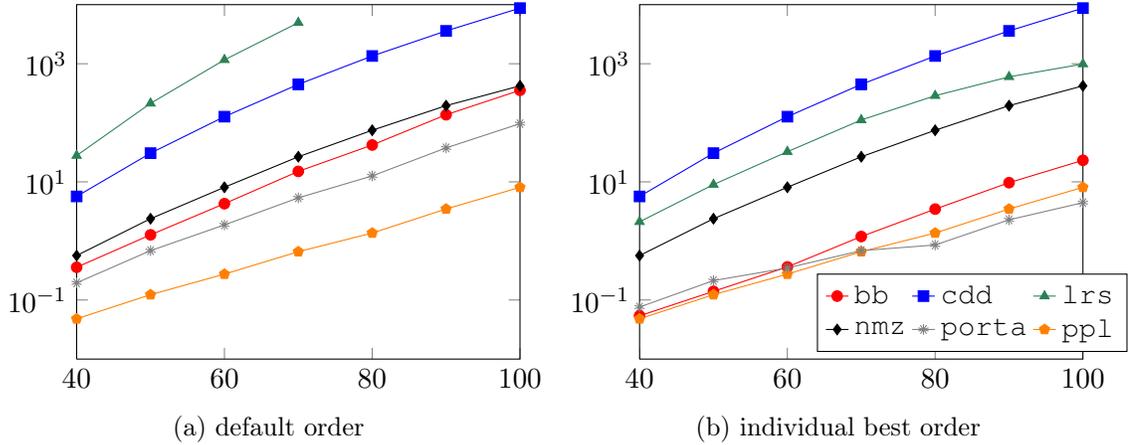

Looking at the timings in Figure~\ref{fig:knapsack:timings:plot:normal} and Table~\ref{tab:knapsack:timings} this experiment seems
to be showing something similar to the previous one: \ppl outperforms the other codes, and \lrs is far behind.
As a major difference to the previous experiment here \porta can nearly compete with \ppl.
For proper comparison, notice that in Figure~\ref{fig:non_sym_cut:timings} the dimension of the input increases with $k$, which is
the parameter on the $x$-axis, while Figure~\ref{fig:knapsack:timings} reports on results in fixed dimension $5$.  The
most interesting aspect of the Fibonacci knapsack examples, however, is that they support the following.

The key to understanding the computational results is the observation that, for $d=6$ and $b=100$, less than $0.5\%$ of
the input points are vertices.  The best possible scenario for any incremental convex hull algorithm is to process a
sequence of points with the vertices coming first.  In more complicated cases the ordering among the vertices may also
matter.
\begin{ruleofthumb}
  For iterative convex hull algorithms the insertion order often makes a difference.  This is particularly true if
  the input contains redundant points.
\end{ruleofthumb}

To illustrate the situation we ran extensive tests on the Fibonacci knapsack examples with various insertion orders; see the rather lengthy Table~\ref{tab:knapsack:timings} in the appendix.
The \cdd and \ppl implementation are almost insensitive to this variation since they employ heuristics to process the input in their preferred order. Also \normaliz orders the input in an initial step, so we did not run tests with varying order for it.
All other codes are very sensitive to the insertion order, even \lrs, which is not an iterative algorithm.
However, the ordering of the input affects the traversal of the search tree.
Figure~\ref{fig:knapsack:timings:plot:best} is an overly condensed version of Table~\ref{tab:knapsack:timings}; for each implementation we took the timing of the best ordering that we found.
\porta and \bb quite clearly want the vertices first, while \lrs prefers a randomized order.
In practice, for a specific input, it will often be unrealistic to tell in advance which ordering will work best for which method.
The purpose of this test is to underline that, in addition to choosing the proper algorithm, it may be worth-while to spend a thought on the input order, too.

\subsection{Voronoi Diagrams}
\label{sec:voronoi}

A completely different convex hull set up comes about as follows.  Consider a configuration $S$ of $m$ points in
$\RR^{d-1}$.  The \emph{Voronoi region} of a \emph{site} $s\in S$ is the set of points in $\RR^{d-1}$ whose distance to
$s$ does not exceed the distance to any other point in $S$.  Conceptually, the distance can be measured in terms of any
fixed metric, but here we are looking at the Euclidean case only.  Then the Voronoi region of $s$ is a convex polyhedron
which may be bounded or not.  An immediate inequality description of a Voronoi region arises from bisecting hyperplanes
between pairs of sites.  Ranging over all sites in $S$ the Voronoi regions form a \emph{polyhedral complex}, \textit{i.e.}, a
set of polyhedra meeting face-to-face, and this polyhedral complex is the \emph{Voronoi diagram} of $S$.  By
construction the Voronoi diagram covers the entire space $\RR^{d-1}$.  Voronoi diagrams are an indispensable tool in
computational geometry and its applications.

The (Euclidean) Voronoi diagram of $S$ can be computed by solving a dual convex hull problem.  This is due to the fact
that the Voronoi diagram of $S$ is the orthogonal projection of the unbounded polyhedron in $\RR^d$ whose $m$ facets are
those tangent hyperplanes to the standard paraboloid $\smallSetOf{(x,\norm{x}^2)}{x\in\RR^{d-1}}$ in $\RR^d$ which lie
above the sites in $S$.  This way the $0$-dimensional cells of the Voronoi diagram of $S$ are the projections of the
vertices of the unbounded polyhedron
\[
V(S) \ := \ \SetOf{(x,\delta)\in\RR^{d-1}\times\RR}{2 \langle s,x\rangle-\norm{s}^2 \le \delta \text{ for all } s\in S} \, ,
\]
which is full-dimensional.  From the vertices and the facets (which are given by the sites) of $V(S)$ one can, for
instance, compute a facet description of each Voronoi region by a combinatorial procedure.  The same is true for other
additional information that might be desired.  In this sense, solving the dual convex hull problem for the polyhedron
$V(S)$ is the key step; e.g., see \cite[\S6.3]{Polyhedral+and+Algebraic+Methods}.

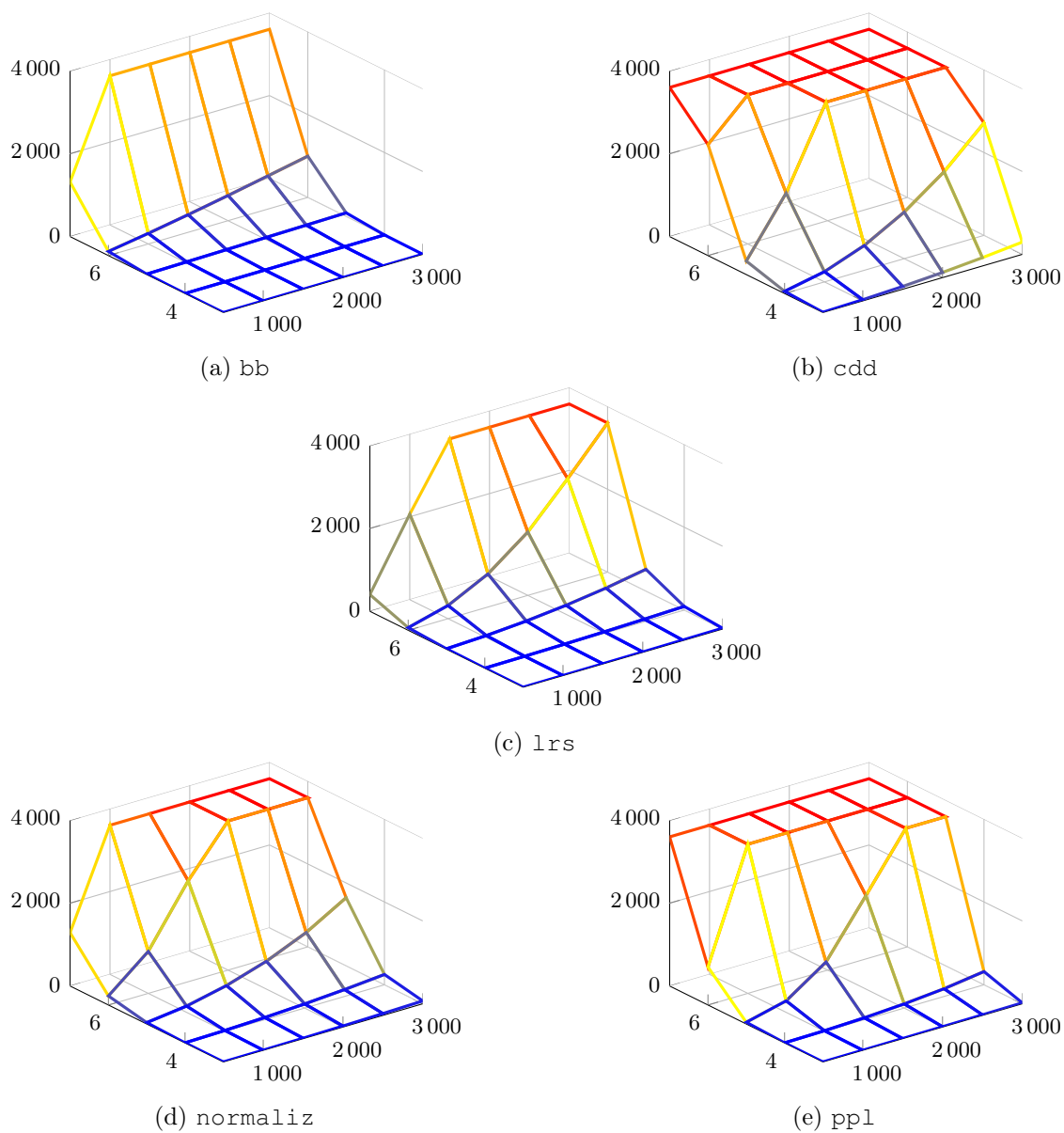
\begin{figure}[p]
  \centering
  \subcaptionbox{\bb}{% This file was created by matlab2tikz v0.4.0.
% Copyright (c) 2008--2013, Nico Schlömer <nico.schloemer@gmail.com>
% All rights reserved.
% 
% The latest updates can be retrieved from
%   http://www.mathworks.com/matlabcentral/fileexchange/22022-matlab2tikz
% where you can also make suggestions and rate matlab2tikz.
% 
% 
% 
\begin{tikzpicture}

% to keep aspect ratio
\setlength{\mywidth}{.33\linewidth} % you might want to change this
\setlength{\myheight}{.85\mywidth} % do not change this

\begin{axis}[%
width=\mywidth,%5.15805555555556in,
height=\myheight,%4.34666666666667in,
view={-37.5}{30},
scale only axis,
xmin=500,
xmax=3000,
xmajorgrids,
ymin=3,
ymax=7,
ymajorgrids,
zmin=0,
zmax=4000,
zmajorgrids,
axis x line*=bottom,
axis y line*=left,
axis z line*=left,
font=\footnotesize,
ticklabel style={/pgf/number format/.cd,1000 sep={\,},}
]
\addplot3[%
mesh,
very thick,
shader=flat,
mesh/rows=6]
table[row sep=crcr,header=false] {
500     3          0.261\\
500     4          1.415\\
500     5         11.144\\
500     6        100.962\\
500     7       1318.488\\
1000    3          0.621\\
1000    4          3.102\\
1000    5         25.280\\
1000    6        255.429\\
1000    7       3600.000\\
1500    3          1.040\\
1500    4          4.875\\
1500    5         40.011\\
1500    6        425.652\\
1500    7       3600.000\\
2000    3          1.517\\
2000    4          6.705\\
2000    5         55.229\\
2000    6        608.085\\
2000    7       3600.000\\
2500    3          2.010\\
2500    4          8.622\\
2500    5         69.848\\
2500    6        800.647\\
2500    7       3600.000\\
3000    3          2.569\\
3000    4         10.547\\
3000    5         85.517\\
3000    6        995.853\\
3000    7       3600.000\\
};
\end{axis}
\end{tikzpicture}%
}
  \hfill
  \subcaptionbox{\cdd}{% This file was created by matlab2tikz v0.4.0.
% Copyright (c) 2008--2013, Nico Schlömer <nico.schloemer@gmail.com>
% All rights reserved.
% 
% The latest updates can be retrieved from
%   http://www.mathworks.com/matlabcentral/fileexchange/22022-matlab2tikz
% where you can also make suggestions and rate matlab2tikz.
% 
% 
% 
\begin{tikzpicture}

% to keep aspect ratio
\setlength{\mywidth}{.33\linewidth} % you might want to change this
\setlength{\myheight}{.85\mywidth} % do not change this

\begin{axis}[%
width=\mywidth,%5.15805555555556in,
height=\myheight,%4.34666666666667in,
view={-37.5}{30},
scale only axis,
xmin=500,
xmax=3000,
xmajorgrids,
ymin=3,
ymax=7,
ymajorgrids,
zmin=0,
zmax=4000,
zmajorgrids,
axis x line*=bottom,
axis y line*=left,
axis z line*=left,
font=\footnotesize,
ticklabel style={/pgf/number format/.cd,1000 sep={\,},}
]
\addplot3[%
mesh,
very thick,
shader=flat,
mesh/rows=6]
table[row sep=crcr,header=false] {
500     3        6.604  \\
500     4       53.436  \\
500     5      326.138  \\
500     6     2679.876  \\
500     7     3600      \\
1000    3       28.420  \\
1000    4      246.447  \\
1000    5     1686.490  \\
1000    6     3600      \\
1000    7     3600      \\
1500    3       67.606  \\
1500    4      599.624  \\
1500    5     3600      \\
1500    6     3600      \\
1500    7     3600      \\
2000    3      124.657  \\
2000    4     1124.073  \\
2000    5     3600      \\
2000    6     3600      \\
2000    7     3600      \\
2500    3      200.765  \\
2500    4     1817.252  \\
2500    5     3600      \\
2500    6     3600      \\
2500    7     3600      \\
3000    3      295.799  \\
3000    4     2712.942  \\
3000    5     3600      \\
3000    6     3600      \\
3000    7     3600      \\
};
\end{axis}
\end{tikzpicture}%
}\\
  \subcaptionbox{\lrs}{% This file was created by matlab2tikz v0.4.0.
% Copyright (c) 2008--2013, Nico Schlömer <nico.schloemer@gmail.com>
% All rights reserved.
% 
% The latest updates can be retrieved from
%   http://www.mathworks.com/matlabcentral/fileexchange/22022-matlab2tikz
% where you can also make suggestions and rate matlab2tikz.
% 
% 
% 
\begin{tikzpicture}

% to keep aspect ratio
\setlength{\mywidth}{.33\linewidth} % you might want to change this
\setlength{\myheight}{.85\mywidth} % do not change this

\begin{axis}[%
width=\mywidth,%5.15805555555556in,
height=\myheight,%4.34666666666667in,
view={-37.5}{30},
scale only axis,
xmin=500,
xmax=3000,
xmajorgrids,
ymin=3,
ymax=7,
ymajorgrids,
zmin=0,
zmax=4000,
zmajorgrids,
axis x line*=bottom,
axis y line*=left,
axis z line*=left,
font=\footnotesize,
ticklabel style={/pgf/number format/.cd,1000 sep={\,},}
]

\addplot3[%
mesh,
very thick,
shader=flat,
mesh/rows=6]
table[row sep=crcr,header=false] {
500     3    0.496      \\
500     4    1.955      \\
500     5    11.744     \\
500     6    63.601     \\
500     7    404.055    \\
1000    3    2.149      \\
1000    4    8.677      \\
1000    5    50.578     \\
1000    6    318.346    \\
1000    7    2063.702   \\
1500    3    5.034      \\
1500    4    20.541     \\
1500    5    119.551    \\
1500    6    798.688    \\
1500    7    3600       \\
2000    3    9.147      \\
2000    4    37.901     \\
2000    5    220.404    \\
2000    6    1530.386   \\
2000    7    3600       \\
2500    3    14.741     \\
2500    4    60.886     \\
2500    5    354.748    \\
2500    6    2537.572   \\
2500    7    3600       \\
3000    3    21.264     \\
3000    4    89.3       \\
3000    5    526.344    \\
3000    6    3600       \\
3000    7    3600       \\
};
\end{axis}
\end{tikzpicture}%
}\\
  \subcaptionbox{\normaliz}{% This file was created by matlab2tikz v0.4.0.
% Copyright (c) 2008--2013, Nico Schlömer <nico.schloemer@gmail.com>
% All rights reserved.
% 
% The latest updates can be retrieved from
%   http://www.mathworks.com/matlabcentral/fileexchange/22022-matlab2tikz
% where you can also make suggestions and rate matlab2tikz.
% 
% 
% 
\begin{tikzpicture}

% to keep aspect ratio
\setlength{\mywidth}{.33\linewidth} % you might want to change this
\setlength{\myheight}{.85\mywidth} % do not change this

\begin{axis}[%
width=\mywidth,%5.15805555555556in,
height=\myheight,%4.34666666666667in,
view={-37.5}{30},
scale only axis,
xmin=500,
xmax=3000,
xmajorgrids,
ymin=3,
ymax=7,
ymajorgrids,
zmin=0,
zmax=4000,
zmajorgrids,
axis x line*=bottom,
axis y line*=left,
axis z line*=left,
font=\footnotesize,
ticklabel style={/pgf/number format/.cd,1000 sep={\,},}
]
\addplot3[%
mesh,
very thick,
shader=flat,
mesh/rows=6]
table[row sep=crcr,header=false] {
500     3     0.605	  \\
500     4     3.108	  \\
500     5     33.374	  \\
500     6     219.029	  \\
500     7     1287.456	  \\
1000    3     3.683	  \\
1000    4     21.365	  \\
1000    5     151.977	  \\
1000    6     1013.544  \\
1000    7     3600	  \\
1500    3     11.179	  \\
1500    4     52.959	  \\
1500    5     365.866	  \\
1500    6     2436.256  \\
1500    7     3600	  \\
2000    3     24.95	  \\
2000    4     100.007	  \\
2000    5     671.886	  \\
2000    6     3600	  \\
2000    7     3600	  \\
2500    3     44.928	  \\
2500    4     165.561	  \\
2500    5     1085.843  \\
2500    6     3600	  \\
2500    7     3600	  \\
3000    3     68.525	  \\
3000    4     249.134	  \\
3000    5     1644.61	  \\
3000    6     3600	  \\
3000    7     3600      \\
};
\end{axis}
\end{tikzpicture}%
}
  \hfill
  \subcaptionbox{\ppl}{% This file was created by matlab2tikz v0.4.0.
% Copyright (c) 2008--2013, Nico Schlömer <nico.schloemer@gmail.com>
% All rights reserved.
% 
% The latest updates can be retrieved from
%   http://www.mathworks.com/matlabcentral/fileexchange/22022-matlab2tikz
% where you can also make suggestions and rate matlab2tikz.
% 
% 
% 
\begin{tikzpicture}

% to keep aspect ratio
\setlength{\mywidth}{.33\linewidth} % you might want to change this
\setlength{\myheight}{.85\mywidth} % do not change this

\begin{axis}[%
width=\mywidth,%5.15805555555556in,
height=\myheight,%4.34666666666667in,
view={-37.5}{30},
scale only axis,
xmin=500,
xmax=3000,
xmajorgrids,
ymin=3,
ymax=7,
ymajorgrids,
zmin=0,
zmax=4000,
zmajorgrids,
axis x line*=bottom,
axis y line*=left,
axis z line*=left,
font=\footnotesize,
ticklabel style={/pgf/number format/.cd,1000 sep={\,},}
]
\addplot3[%
mesh,
very thick,
shader=flat,
mesh/rows=6]
table[row sep=crcr,header=false] {
500     3     0.289     \\
500     4     2.508     \\
500     5     41.547    \\
500     6     917.853   \\
500     7     3600      \\
1000    3     1.355     \\
1000    4     14.48     \\
1000    5     277.779   \\
1000    6     3600      \\
1000    7     3600      \\
1500    3     3.757     \\
1500    4     45.456    \\
1500    5     935.003   \\
1500    6     3600      \\
1500    7     3600      \\
2000    3     7.876     \\
2000    4     100.803   \\
2000    5     2255.994  \\
2000    6     3600      \\
2000    7     3600      \\
2500    3     13.947    \\
2500    4     186.21    \\
2500    5     3600      \\
2500    6     3600      \\
2500    7     3600      \\
3000    3     22.495    \\
3000    4     319.268   \\
3000    5     3600      \\
3000    6     3600      \\
3000    7     3600      \\
};
\end{axis}
\end{tikzpicture}%
}
  \caption{Running times (in seconds) for Voronoi diagrams of random point sets depending on the dimension $d \in [3,7]$ and the number of points $m \in [500,3000]$. See Table~\ref{tab:voronoi:timings} for exact timings.}
  \label{fig:voronoi:timings}
\end{figure}

For our experiment we choose $m$ points uniformly at random within the cube $[-1,1]^{d-1}$.
The polytopes $V(S)$ were given to \polymake in reduced facet description via \texttt{FACETS}, derived from the sites given by the random points. Except for \bb where we added an explicit polarization step because this implementation originally worked for primal convex hull only, thus the input was \texttt{VERTICES}.  See the accompanying web page~\cite{polymake_exp} for full details.
For each choice of
parameters we do this experiment ten times.  The timings in Table~\ref{tab:voronoi:timings} give the average values.  We
tried $m\in\{500,1000,1500,\dots,3000\}$ and $3 \le d \le 7$, but we restricted our attention to those cases which could
be computed within one hour of single-threaded CPU time.  While Table~\ref{tab:voronoi:timings} also lists a few timings
which are longer, everything is truncated to one hour in Figure~\ref{fig:voronoi:timings}.  Notice that the polyhedron
$V(S)$ is simple (\textit{i.e.}, each vertex is contained in $d$ facets) almost surely.

On this non-degenerate input \cdd does not work so well and can only solve a few small problems.  Reverse search is a
method which is designed for non-degenerate input, and so its implementation \lrs performs quite well.  In fact, the
worst case running time of reverse search is bounded by $O(dmn)$ on non-degenerate input; and this fits with our experiment.
It is plausible (but not proved, as far as we know) that $n$ grows only slightly larger than linearly in $m$ in this
case, and so the running time in fixed dimension approximately grows quadratically.  Yet this is outperformed by the
beneath-and-beyond method \bb, implemented in \polymake.  Its worst case running time was estimated at $O(d^5mt^2)$ in
\cite{Beneath-Beyond}, where $t$ is the size of a triangulation of the polar of $V(S)$ which \bb implicitly computes
on the way.  While it is plausible (but again not proved, as far as we know) that $t$ is small in this case, possibly
even about linear in $m$, the above bound seems to be too pessimistic.  The empirical running times on this input grow
almost linearly in fixed dimension.  Notice also that, for any $d$ and very small $m$ \lrs, \normaliz and \ppl beat \bb.

\begin{ruleofthumb}
  On random input the beneath-and-beyond algorithm often behaves very well.
\end{ruleofthumb}

% \begin{itemize}
%   \item $d=5$ and $\text{Voronoi-sites}=m=500$; Vertices $n$ are calculated
%   \item Median: 12816
%   \item Expected Value: 12816.352
%   \item Var: 8655.964
%   \item Standard Deviation: 93.037
%   \item Max: 13130
%   \item Min: 12538
% \end{itemize}

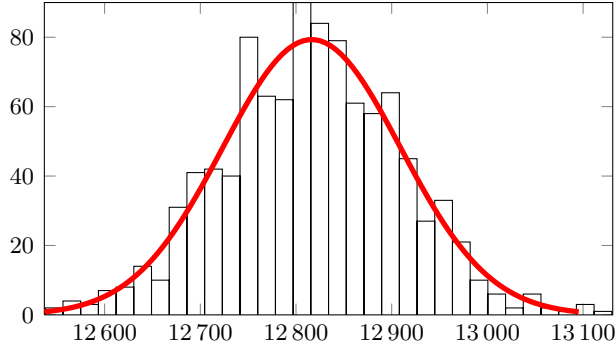
\begin{figure}[ht]
  \centering
  % This file was created by matlab2tikz v0.4.5 running on MATLAB 8.2.
% Copyright (c) 2008--2014, Nico Schlömer <nico.schloemer@gmail.com>
% All rights reserved.
% 
% The latest updates can be retrieved from
%   http://www.mathworks.com/matlabcentral/fileexchange/22022-matlab2tikz
% where you can also make suggestions and rate matlab2tikz.
% 
\begin{tikzpicture}

% to keep aspect ratio
\setlength{\mywidth}{.5\linewidth} % you might want to change this
\setlength{\myheight}{.55\mywidth} % do not change this

\begin{axis}[%
width=\mywidth, %6.44972222222222in,
height=\myheight, %5.08016666666667in,
scale only axis,
xmin=12537,
xmax=13131,
ymin=0,
ymax=90,
font=\footnotesize,
scaled x ticks = false,
xticklabel style={/pgf/number format/.cd,1000 sep={\,},}
]

\addplot[area legend,forget plot]
 table[row sep=crcr, point meta=\thisrow{c}] {
x	y	c\\
12538	0	1	\\
12538	2	1	\\
12556.5	2	1	\\
12556.5	0	1	\\
12556.5	0	1	\\
12556.5	4	1	\\
12575	4	1	\\
12575	0	1	\\
12575	0	1	\\
12575	3	1	\\
12593.5	3	1	\\
12593.5	0	1	\\
12593.5	0	1	\\
12593.5	7	1	\\
12612	7	1	\\
12612	0	1	\\
12612	0	1	\\
12612	8	1	\\
12630.5	8	1	\\
12630.5	0	1	\\
12630.5	0	1	\\
12630.5	14	1	\\
12649	14	1	\\
12649	0	1	\\
12649	0	1	\\
12649	10	1	\\
12667.5	10	1	\\
12667.5	0	1	\\
12667.5	0	1	\\
12667.5	31	1	\\
12686	31	1	\\
12686	0	1	\\
12686	0	1	\\
12686	41	1	\\
12704.5	41	1	\\
12704.5	0	1	\\
12704.5	0	1	\\
12704.5	42	1	\\
12723	42	1	\\
12723	0	1	\\
12723	0	1	\\
12723	40	1	\\
12741.5	40	1	\\
12741.5	0	1	\\
12741.5	0	1	\\
12741.5	80	1	\\
12760	80	1	\\
12760	0	1	\\
12760	0	1	\\
12760	63	1	\\
12778.5	63	1	\\
12778.5	0	1	\\
12778.5	0	1	\\
12778.5	62	1	\\
12797	62	1	\\
12797	0	1	\\
12797	0	1	\\
12797	90	1	\\
12815.5	90	1	\\
12815.5	0	1	\\
12815.5	0	1	\\
12815.5	84	1	\\
12834	84	1	\\
12834	0	1	\\
12834	0	1	\\
12834	79	1	\\
12852.5	79	1	\\
12852.5	0	1	\\
12852.5	0	1	\\
12852.5	61	1	\\
12871	61	1	\\
12871	0	1	\\
12871	0	1	\\
12871	58	1	\\
12889.5	58	1	\\
12889.5	0	1	\\
12889.5	0	1	\\
12889.5	64	1	\\
12908	64	1	\\
12908	0	1	\\
12908	0	1	\\
12908	45	1	\\
12926.5	45	1	\\
12926.5	0	1	\\
12926.5	0	1	\\
12926.5	27	1	\\
12945	27	1	\\
12945	0	1	\\
12945	0	1	\\
12945	33	1	\\
12963.5	33	1	\\
12963.5	0	1	\\
12963.5	0	1	\\
12963.5	21	1	\\
12982	21	1	\\
12982	0	1	\\
12982	0	1	\\
12982	10	1	\\
13000.5	10	1	\\
13000.5	0	1	\\
13000.5	0	1	\\
13000.5	6	1	\\
13019	6	1	\\
13019	0	1	\\
13019	0	1	\\
13019	2	1	\\
13037.5	2	1	\\
13037.5	0	1	\\
13037.5	0	1	\\
13037.5	6	1	\\
13056	6	1	\\
13056	0	1	\\
13056	0	1	\\
13056	2	1	\\
13074.5	2	1	\\
13074.5	0	1	\\
13074.5	0	1	\\
13074.5	1	1	\\
13093	1	1	\\
13093	0	1	\\
13093	0	1	\\
13093	3	1	\\
13111.5	3	1	\\
13111.5	0	1	\\
13111.5	0	1	\\
13111.5	1	1	\\
13130	1	1	\\
13130	0	1	\\
};

\addplot [color=red,solid,line width=2.0pt,forget plot]
  table[row sep=crcr]{
12537.2397403737	0.881250714356249	\\
12542.8783506771	1.05502851719306	\\
12548.5169609805	1.25844350334142	\\
12554.1555712838	1.49557451217514	\\
12559.7941815872	1.77087209770559	\\
12565.4327918906	2.08915730987335	\\
12571.071402194	2.45561285522585	\\
12576.7100124974	2.87576544836988	\\
12582.3486228007	3.35545819522885	\\
12587.9872331041	3.90081192092283	\\
12593.6258434075	4.51817447392948	\\
12599.2644537109	5.21405720794806	\\
12604.9030640143	5.99505806614951	\\
12610.5416743176	6.86777097007986	\\
12616.180284621	7.83868154626053	\\
12621.8188949244	8.91404960409325	\\
12627.4575052278	10.0997792030947	\\
12633.0961155312	11.4012776073244	\\
12638.7347258345	12.8233049089605	\\
12644.3733361379	14.3698165977429	\\
12650.0119464413	16.0438018425433	\\
12655.6505567447	17.8471207177919	\\
12661.2891670481	19.7803440317751	\\
12666.9277773514	21.8425997758864	\\
12672.5663876548	24.031430493878	\\
12678.2049979582	26.3426660489482	\\
12683.8436082616	28.7703163266953	\\
12689.482218565	31.3064883390623	\\
12695.1208288683	33.9413319771726	\\
12700.7594391717	36.6630182929392	\\
12706.3980494751	39.4577536688589	\\
12712.0366597785	42.309832566738	\\
12717.6752700818	45.2017307394124	\\
12723.3138803852	48.1142398609912	\\
12728.9524906886	51.0266435031363	\\
12734.591100992	53.9169332848678	\\
12740.2297112954	56.762062883826	\\
12745.8683215987	59.5382364538946	\\
12751.5069319021	62.2212268859793	\\
12757.1455422055	64.7867183152505	\\
12762.7841525089	67.21066635806	\\
12768.4227628123	69.4696687919025	\\
12774.0613731156	71.5413388051289	\\
12779.699983419	73.4046725671102	\\
12785.3385937224	75.040402725216	\\
12790.9772040258	76.4313295351794	\\
12796.6158143292	77.5626216806654	\\
12802.2544246325	78.4220794314248	\\
12807.8930349359	79.0003536133805	\\
12813.5316452393	79.2911148957948	\\
12819.1702555427	79.2911691094715	\\
12824.8088658461	79.0005156581149	\\
12830.4474761494	78.4223475292484	\\
12836.0860864528	77.5629929044101	\\
12841.7246967562	76.4317998616779	\\
12847.3633070596	75.0409671078765	\\
12853.001917363	73.4053250260351	\\
12858.6405276663	71.5420725325285	\\
12864.2791379697	69.4704762702055	\\
12869.9177482731	67.2115394878601	\\
12875.5563585765	64.7876485500506	\\
12881.1949688798	62.2222053706115	\\
12886.8335791832	59.5392541635033	\\
12892.4721894866	56.7631107603995	\\
12898.11079979	53.9180023684833	\\
12903.7494100934	51.0277250553388	\\
12909.3880203967	48.1153254780432	\\
12915.0266307001	45.2028124535173	\\
12920.6652410035	42.3109029337405	\\
12926.3038513069	39.4588058413488	\\
12931.9424616103	36.6640460781105	\\
12937.5810719136	33.94232987921	\\
12943.219682217	31.3074515862246	\\
12948.8582925204	28.7712408836312	\\
12954.4969028238	26.343548615104	\\
12960.1355131272	24.0322684891328	\\
12965.7741234305	21.8433913147713	\\
12971.4127337339	19.7810878877987	\\
12977.0513440373	17.8478162795031	\\
12982.6899543407	16.0444490630196	\\
12988.3285646441	14.3704159391801	\\
12993.9671749474	12.8238572840183	\\
12999.6057852508	11.4017843187715	\\
13005.2443955542	10.1002418831915	\\
13010.8830058576	8.91447015499608	\\
13016.521616161	7.83906208256765	\\
13022.1602264643	6.86811376445741	\\
13027.7988367677	5.9953654987938	\\
13033.4374470711	5.21433172045583	\\
13039.0760573745	4.5184185279722	\\
13044.7146676779	3.90102796195529	\\
13050.3532779812	3.35564862131256	\\
13055.9918882846	2.87593258405854	\\
13061.630498588	2.45575893034721	\\
13067.2691088914	2.0892844429483	\\
13072.9077191947	1.7709822836723	\\
13078.5463294981	1.49566961403265	\\
13084.1849398015	1.25852524729167	\\
13089.8235501049	1.05509849084212	\\
13095.4621604083	0.881310367531189	\\
};
\end{axis}
\end{tikzpicture}%
  \caption{Histogram of the numbers of vertices in 1000 random Voronoi diagram experiments, the boundaries of the $x$-axis correspond to the minimum and maximum values.}
  \label{fig:voronoi:stats:n}
\end{figure}

% beneath and beyond (5,500)
% \begin{itemize}
%   \item Median: 11.29
%   \item Expected Value: 11.29
%   \item Var: 0.023588008
%   \item Standard Deviation: 0.15358387
%   \item Max: 11.8
%   \item Min: 10.83
% \end{itemize}
%
% cdd (5,500)
% \begin{itemize}
%   \item Median: 332.81
%   \item Expected Value: 332.84245
%   \item Var: 11.5477286261
%   \item Standard Deviation: 3.398194907
%   \item Max: 343.43
%   \item Min: 320.84
% \end{itemize}
%
% lrs (5,500)
% \begin{itemize}
%   \item Median: 11.8
%   \item Expected Value: 11.80157
%   \item Var: 0.0338022373
%   \item Standard Deviation: 0.1838538478
%   \item Max: 12.39
%   \item Min: 11.3
% \end{itemize}
%
% ppl (5,500)
% \begin{itemize}
%   \item Median: 47.63
%   \item Expected Value: 47.64231
%   \item Var: 0.6570762401
%   \item Standard Deviation: 0.8106023934
%   \item Max: 50.32
%   \item Min: 45.2500000001
% \end{itemize}
%

% REMOVED since 1st reviewer did not like it.
% temporarily creates unresolved references

% \begin{figure}[p]
%   \subcaptionbox{\bb}{\input{img/histfit_bb.tikz}}\hfill
%   \subcaptionbox{\cdd}{\input{img/histfit_cdd.tikz}}\\
%   \subcaptionbox{\lrs}{\input{img/histfit_lrs.tikz}}\hfill
%   \subcaptionbox{\ppl}{\input{img/histfit_ppl.tikz}}

%   \caption{Histograms of timings of $1000$ random Voronoi diagram experiments, the boundaries of the $x$-axis correspond to the minimum and maximum values.}
%   \label{fig:voronoi:stats:time}
% \end{figure}

We conclude the report on this experiment with an explanation why it suffices to list the average over only ten runs per
parameter set via a simple statistical analysis.  To this end we fix one set of parameters for which we do the random sampling
one thousand times.  We pick $d=5$ and $m=500$ because all algorithms terminate within a reasonable time-span.

The first random variable to look at is the number $n$ of vertices of the Voronoi diagram, whose distribution is plotted
in Figure~\ref{fig:voronoi:stats:n}.  The numbers range from 12\thinspace538 to 13\thinspace130 with average value
12\thinspace816.352, median 12\thinspace816 and standard deviation 93.037.  For comparison the diagram also shows a
normal distribution with the given expected value and standard deviation.  To get an idea, it follows from the
Generalized Lower Bound Theorem of Murai and Nevo~\cite{MuraiNevo13} that $n \ge 1\thinspace982$ (provided that the
input is in general position).  On the other hand McMullen's Upper Bound Theorem \cite{McMullen71} yields $n \le
246\thinspace512$.

The other natural random variables are the running times for the various convex hull codes.
Like for the output sizes, these are strongly concentrated around their respective mean.
We omit the details.

\begin{remark}
  We do not report \porta timings here since \porta dies with an error message on each input of this kind.
  This seems to be related to a flaw in \porta's implementation of the arithmetic.
  % "invalid denominator"
  The most recent version \cite{porta} of \porta was published in 2009.
\end{remark}

\subsection{Symmetric Cut Polytopes}
\label{sec:symmetric-cut}

We continue to examine primal convex hull computations arising from cut polytopes.  The difference to
Section~\ref{sec:non-symmetric-cut} is that here we look into graphs with non-trivial symmetry.  A symmetric convex hull
problem can be solved in various ways by graph based search techniques up to symmetry which need to be combined with
some standard convex hull algorithm applied to subproblems of smaller size; see Bremner \textit{et al.} \cite{Bremner09}.
This has been implemented in \sympol by Rehn \cite{RehnDipl}, and this is also available in \polymake, through which \sympol can
currently be combined with \bb, \cdd, \lrs and \ppl.

Our first example is the path $P_k$ on $k$ nodes, which has $k-1$ edges.  The cut polytope
$\Cut(P_k)$ is the $0/1$-cube of dimension $k-1$.  It has $n=2^{k-1}$ vertices and $m=2k-2$ facets.  The graph $P_k$
is symmetric with respect to exchanging both ends, which yields an automorphism group of order two.  Clearly, the
$0/1$-cube has a much larger group of automorphisms (of order $2^{k-1}\cdot(k-1)!$), but most of its elements are not induced by graph automorphisms.

The second example is the cycle $C_k$ of length $k$.  The cut polytope $\Cut(C_k)$ is the convex hull of all
$0/1$-vectors of length $k$ with an even number of $1$'s.  It has $n=2^{k-1}$ vertices and $m=2k+2^{k-1}$ facets.  The
automorphism group of the graph $C_k$ is the dihedral group of order $2k$.

The third example is the complete graph $K_k$ on $k$ nodes, which has $\tbinom{k}{2}=k(k-1)/2$ edges.  The number of
vertices of $\Cut(K_k)$ equals $2^{k-1}$.  The facet descriptions are known for $k\le 8$; see
\cite{Bejraburnin,smapo,MR2361926}. For $k=6$ the cut polytope has $368$ facets, for $k=7$ it has $116\thinspace764$
facets, for $k=8$ it has $217\thinspace093\thinspace472$ facets, and at least
$12\thinspace246\thinspace651\thinspace158\thinspace320$ facets for $k=9$~\cite{Bejraburnin}.  The latter number is
conjectured to be the true number of facets.  The automorphism group of $K_k$ is the full symmetric group of degree $k$.

As for the case of non-symmetric cut polytopes the input was given as a reduced list of vertices. However, in this case
we also passed a list of generators of the symmetry group to \polymake for the computations. See~\cite{polymake_exp} for the complete setup.

 {\footnotesize
\begin{table}[t]
  \centering
  \caption{Computation times (in seconds) for the convex hull of several cut polytopes where the symbol \textoom stands for \emph{out of memory} and $*$ means that we have only done one iteration of the test.}
  \label{tab:cutpolytopes}
  \begin{tabular}[t]{l@{\hspace{.1cm}}cccd{2}d{2}d{2}d{2}d{2}d{2}d{2}d{2}}
    \toprule
    $G$       & $d$ & $n$ & $m$ & \multicolumn{1}{c}{\text{\bb}}      & \multicolumn{1}{c}{\text{\cdd}}      & \multicolumn{1}{c}{\text{\lrs}} & \multicolumn{1}{c}{\text{\ppl}} & \multicolumn{1}{c}{\text{\symbb}} & \multicolumn{1}{c}{\text{\symcdd}} & \multicolumn{1}{c}{\text{\symlrs}} & \multicolumn{1}{c}{\text{\symppl}}\\
    \midrule
    $P_9$    &  8 & 256 &     16 &  0.38 &     0.18   &      5.36   & 0.07    & 0.29 & 0.23 &  1.36  & 0.08 \\
    $P_{10}$  &  9 & 512 &     18 &  3.09 &     0.52  &     90.80   & 0.07 &    2.04 & 0.95 & 26.12   & 0.25 \\[1.2ex]
    $C_9$    &  9 & 256 &    274 &  2.05 &     0.88   &     29.46    & 0.09    & 1.07 & 0.97 &  2.86  & 0.56 \\
    $C_{10}$  & 10 & 512 &    532 & 24.99 &     4.70  &    602.04   & 0.20     & 6.76 & 4.22 & 62.70  &  2.40\\[1.2ex]
    $K_6$    & 15 &  32 &    368 &  6.80 &     0.25   &      1.74   & 0.09     & 1.04 & 0.20 &  0.20       & 0.10 \\
    $K_7$    & 21 &  64 & 116764 &  \oom  & 14329.75^* &  31309.01^* & 2212.73^* & \oom  & 129.98 & 1217.96 &  32.00\\
%              d     n       m       bb        cdd           lrs            ppl       symbb   symcdd   symlrs     symppl
    \bottomrule
  \end{tabular}\bigskip
\end{table}
}

Table~\ref{tab:cutpolytopes} gives the average time in seconds used by the algorithms in \bb, \cdd, \lrs, \ppl and their
combinations with \sympol for the convex hull computations of the three classes of cut polytopes.  Independent of the
algorithm, at least for large input the advantage of taking symmetry into account is evident.  The larger the group of
automorphisms the greater is the gain.

\begin{ruleofthumb}
  If the input is large but symmetric do compute convex hulls up to symmetry.
\end{ruleofthumb}

The differences in performance of the convex hulls codes on this input are consistent with the non-symmetric cut polytopes from Section~\ref{sec:non-symmetric-cut}.
We take this as an indication that the rule of thumb above is independent of the choice of the basic algorithm.

\section{Integer Points}
\label{sec:integer_points}

\noindent
In this section we compare various methods to count the integer points in a polytope. This is a fundamental task which
appears in various different areas of mathematics, among them number theory, statistics, algebraic geometry, and
representation theory (see \textit{e.g.}~\cite{DeLoera05} for an overview).  In the context of optimization this most
frequently occurs as a first step if one is interested in computing the facets of integer hulls.

\subsection{Complexity status}
While this was already mentioned in the introduction, we wish to repeat: Even to decide if a polyhedron contains some
integer point is NP-complete \cite{GareyJohnson79}.  On top of this there is no polynomial bound on the number of
lattice points with respect to the input size.  Hence all algorithms that enumerate the lattice points easily show an
(at least) exponential worst case behavior with respect to time and memory. Also note that in dimensions greater than three, a
polytope without any lattice points, except for the vertices, can have arbitrarily large volume.
In contrast, for a fixed positive number of interior lattice points the volume is bounded.

Despite these obstacles the number of lattice points can be determined in polynomial time in fixed dimension with Barvinok's
algorithm~\cite{barvinok_algo}.
For a short explanation of the algorithm see the next section.
Variations of this algorithm have been implemented in two different software packages.
This leads to the first rule of thumb for this section.

\begin{ruleofthumb}\label{rot:lattice:latte}
  If you do not know anything about your input and you are only interested in the number of lattice points, try
  {\rm \latte}~\cite{latte} or {\rm \barvinok}~\cite{SkimoBarvinok}.
\end{ruleofthumb}

The reason for this general recommendation is that in practice any method which explicitly enumerates the lattice points will often feel the limit of the amount of memory available.

\subsection{Common Algorithms}

Algorithms to solve the enumeration or counting problem in practice basically fall into four different categories.

The basic form of Barvinok's algorithm computes the rational multivariate generation function for the lattice points in
an affine polyhedral cone $C$ using a \emph{signed} decomposition of $C$ into simplicial unimodular cones whose
generating functions are easy to determine. The key feature of this algorithm is that we can bound the number of
unimodular cones and their faces in this decomposition by a polynomial in the input size if the dimension is fixed. In
its original form the algorithm has to care for lattice points in intersections of cones in the decomposition using
inclusion-exclusion. There are two commonly used variations that avoid this, either working with the dual of the cone or
using \emph{irrational decomposition}~\cite{koeppeirrational}. The first approach is based on the observation that duals
of low-dimensional cones are non-pointed and the rational generating function of a non-pointed cone vanishes. In the
second approach one translates the cone by a small vector in such a way that no lattice point is contained in a
lower-dimensional face (this can be done without explicitly computing the shift). See~\cite{barvinok_algo} for details.

By Brion's Theorem the generating function of the integer points in a rational polytope $P$ is given by the the sum of the
generating functions of all vertex cones. Via Barvinok's algorithm we obtain an algorithm to compute a polynomial size
representation of the generating function of $P$ in polynomial time, if the dimension is fixed. Evaluation at $1$
returns the number of lattice points. This algorithm has been implemented in the software packages {\rm
  \latte}~\cite{DeLoeraEtAl13,latte} and {\rm \barvinok}~\cite{SkimoBarvinok}. The latter package also contains several
modifications of the original algorithm (see, \textit{e.g.}, \cite{KoeppeVerdoolaegeWoods2008}).
Here we only test \latte, since currently there is no interface between \polymake and \barvinok.
% There basically is, but it ends in a segfault I currently cannot track down.

Historically, \latte uses decomposition in dual space as the default method for triangulating the vertex cones. However,
the current version \lattem also implements irrational decomposition and a mix of dual and irrational decomposition. It also
allows to stop the signed decomposition if the determinant of a cone falls below a threshold. This is motivated by the
observation that determinants initially drop fast in the decomposition process, but breaking down cones with small
determinant into unimodular ones often takes a considerable amount of time.
\latte also allows to compute the generating function of the homogenization cone $\sigma_P:=\cone(P\times \{1\})$ via Barvinok's algorithm.

By default, \polymake calls the default version of \latte using decomposition in dual space. This can be changed by
setting command line options for the call to \latte. In \polymake, this can be done with the command
\texttt{set\_custom}, e.g.\
\begin{lstlisting}
polytope> set_custom($latte_count_param="--irrational-all-primal --maxdet=25 --exponential");
\end{lstlisting}%$
For more options see the manual of \lattem.
The parameters chosen may have a tremendous influence on the running time.
As for the convex hulls there is no a priori way to determine which algorithm works best on a given instance.
The manual of \lattem gives some hints.
We demonstrate this in our experiments by running them once with the default options and once with the above choice (denoted by \latte{}$^\vee$).

The second type of algorithm directly enumerates the lattice points of a polytope. Two such methods are implemented in
\polymake.
\begin{enumerate}
\item The algorithm \bbox encloses a polytope into a scaled cube with edges parallel to the
  coordinate axes. It then checks for every point in the cube whether it is contained in the polytope.
\item The algorithm \projection recursively enumerates lattice points in the fibers over
  lattice points in the projection of the polytope into one of the coordinate hyperplanes. This uses an optimized
  Fourier--Motzkin elimination to create irredundant descriptions for the polytopes which are the images of the
  projection.
\end{enumerate}

The third method that can be used to enumerate lattice points in fact solves the more general problem of finding a
Hilbert basis in the cone $\sigma_P$. As all lattice points in $P\times \{1\}$ must be in the Hilbert basis we can
filter out all generators with a higher last coordinate to obtain a list of lattice points. Various algorithms exist for
enumeration of Hilbert bases, but there are only two large software projects that have implemented some of the
methods. To both software projects \polymake provides an interface.
\begin{enumerate}
\item The program \ftit\cite{4ti2} uses a project-and-lift approach to compute Hilbert bases of polyhedral cones
  intersected with arbitrary lattices~\cite{4ti2:hilbert}. Computation of integer points in polyhedra is treated as a
  special case by bounding the height of the generators.
\item The program \normaliz~\cite{normaliz2} computes (among other things) Hilbert bases of cones and enumerates lattice
  points in polytopes. Roughly speaking, the algorithm constructs a triangulation of a cone (\textit{e.g.}, via a
  placing triangulation similar to the beneath-and-beyond method explained above), then computes a Hilbert basis of each
  simplicial cone and finally reduces the union of all Hilbert bases of the simplicial cones to a Hilbert basis of the
  original cone. Computation of a Hilbert basis of a simplicial cone is done by enumerating all lattice points in the
  fundamental parallelepiped using a transformation to the positive orthant. For the enumeration of the lattice points
  of a polytope \normaliz considers the homogenization cone and discards all generators of height 2 or above in the
  computation of the Hilbert basis.

  There exists an option to run \normaliz in a \emph{dual mode} which is based on an algorithm by Pottier~\cite{nmzdual}.
  As for the convex hulls there is no a priori way to determine which algorithm works best on a given instance.
  Still, the authors recommend this mode when only an outer description is given as input.
  We demonstrate this in our experiments by running them once with the default options and once in dual mode (denoted by \normaliz{}$^\vee$).
\end{enumerate}

Generally, as Barvinok's algorithm is polynomial in fixed dimension it should be used if only the number of lattice
points needs to be determined. We will see an example supporting this recommendation in the following section, where we
pick up the fractional knapsack polytopes from Section~\ref{sec:knapsack:hull} again and discuss the computation of the lattice
points used as input above. However, in the subsequent section we will introduce \emph{random box} polytopes and show that \latte fails on those, seemingly simply structured, polytopes.
At the same time the theoretically inferior methods work quite well.
We will shortly discuss further examples in Section~\ref{sec:lattice:additional_comments}.

\medskip

Clearly, one can tailor implementations to address special cases.
For purposes related to integer linear programming it is most interesting to count or to enumerate the $0/1$-points in a polytope which is given in terms of inequalities.
The prototypical implementation is \azove \cite{azove}, which is also interfaced in \polymake (via text file exchange).
\azove traverses all paths in a binary decision diagram to enumerate the $0/1$-points.
\begin{remark}\label{rem:azove}
  The latest version of \azove was published in 2007.
  To use it with a modern Linux kernel requires some patching, as otherwise the initial memory allocation fails.
\end{remark}
% It is the only software in our tests where we used a file based interface, because there is no library version of \azove.
%Thus, in our timings we have a slight overhead due to writing and reading the input/output files.

\subsection{Knapsack Polytopes}
\label{sec:knapsack:integer}
We evaluate the different algorithms for counting and enumerating the lattice points in Fibonacci fractional knapsack polytopes as
introduced in Equation~\eqref{eq:knapsack} of Section~\ref{sec:session}.  This enumeration is the first step of the
integer hull computation in Section~\ref{sec:knapsack:hull}.  According to our Rule of Thumb~\ref{rot:lattice:latte} we
would expect \latte to perform significantly better than the other algorithms, but as we will see, the picture is
slightly more complex.
\begin{remark}
  By adding a slack variable a fractional knapsack polytope can also be written as the set of non-negative solutions to a single
  linear equation. In this context Baldoni \textit{et al.}~\cite{BaldoniBerlineDeLoeraDutraKoeppeVergne:1312.7147}
  described a polynomial time algorithm to compute the top coefficients of a quasi-polynomial which counts the integer
  solutions.
\end{remark}

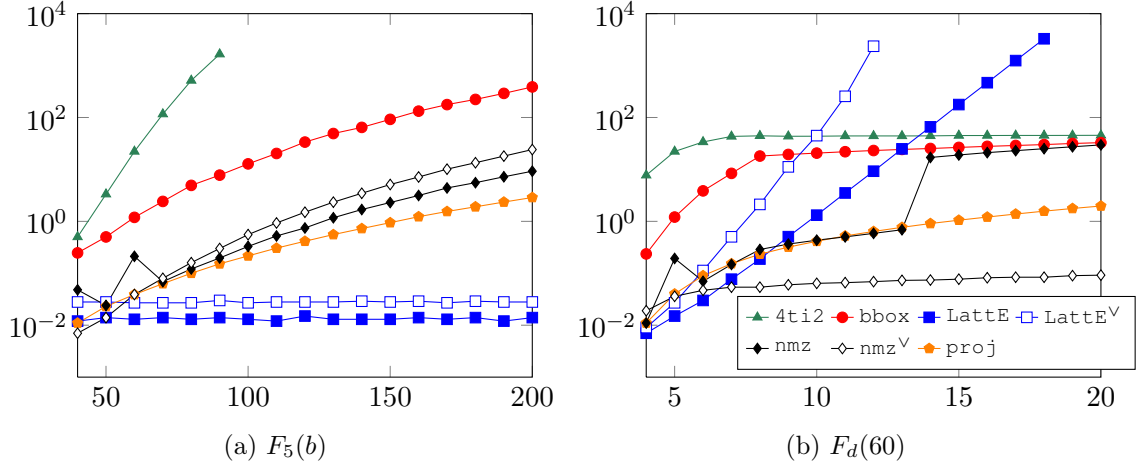
\begin{figure}[tb]
  \centering
  \subcaptionbox{$F_5(b)$\label{fig:knapsack:lattice:timings:plot:rhs}}{% This file was created by matlab2tikz v0.4.5 running on MATLAB 8.2.
% Copyright (c) 2008--2014, Nico Schlömer <nico.schloemer@gmail.com>
% All rights reserved.
% 
% The latest updates can be retrieved from
%   http://www.mathworks.com/matlabcentral/fileexchange/22022-matlab2tikz
% where you can also make suggestions and rate matlab2tikz.
%
%
% defining custom colors
\definecolor{mycolor1}{rgb}{0.16863,0.50588,0.33725}%
\begin{tikzpicture}

% to keep aspect ratio
\setlength{\mywidth}{.4\linewidth} % you might want to change this
\setlength{\myheight}{.8\mywidth} % do not change this

\begin{axis}[%
width=\mywidth, %6.4cm,
height=\myheight, %5cm,
scale only axis,
xmin=40,
xmax=200,
ymode=log,
ymin=0.001,
ymax=10000,
yminorticks=true,
legend style={empty legend}%{at={(0.74,0.02)},anchor=south west,draw=black,fill=white,legend cell align=left}
]
\addplot [ftitstyle]
  table[row sep=crcr]{
40  0.500   \\   
50  3.331   \\   
60  22.171  \\  
70  116.258 \\ 
80  515.834 \\ 
90  1657.244\\
};
%\addlegendentry{\ftit};

\addplot [bboxstyle]
  table[row sep=crcr]{
40     0.247 \\  
50     0.501 \\  
60     1.193 \\  
70     2.420 \\  
80     4.928 \\  
90     7.766 \\  
100   12.787 \\  
110   20.288 \\ 
120   33.606 \\ 
130   49.113 \\ 
140   64.496 \\ 
150   92.070 \\ 
160  132.536 \\ 
170  177.143 \\
180  222.824 \\
190  291.333 \\
200  387.886 \\
};	
%\addlegendentry{\bbox};

\addplot [lattestyle]
  table[row sep=crcr]{
40   0.012 \\
50   0.014 \\
60   0.013 \\
70   0.014 \\
80   0.013 \\
90   0.014 \\
100  0.013 \\
110  0.012 \\
120  0.015 \\
130  0.013 \\
140  0.013 \\
150  0.013 \\
160  0.014 \\
170  0.013 \\
180  0.014 \\
190  0.012 \\
200  0.014 \\
};
%\addlegendentry{\latte};

\addplot [lattestylealt]
  table[row sep=crcr]{
40   0.028 \\
50   0.028 \\
60   0.027 \\
70   0.027 \\
80   0.027 \\
90   0.030 \\
100  0.027 \\
110  0.028 \\
120  0.028 \\
130  0.028 \\
140  0.029 \\
150  0.028 \\
160  0.029 \\
170  0.027 \\
180  0.029 \\
190  0.028 \\
200  0.028 \\
};
%\addlegendentry{\latte *};

\addplot [projstyle]
  table[row sep=crcr]{
40    0.011\\
50    0.023\\
60    0.039\\
70    0.063\\
80    0.100\\
90    0.152\\
100   0.216\\
110   0.305\\
120   0.414\\
130   0.556\\
140   0.727\\
150   0.943\\
160   1.229\\
170   1.549\\
180   1.907\\
190   2.352\\
200   2.880\\
};
%\addlegendentry{\projection};

\addplot [nmzstyle]
  table[row sep=crcr]{
40   0.048 \\
50   0.024 \\
60   0.212 \\
70   0.069 \\
80   0.121 \\
90   0.197 \\
100  0.329 \\
110  0.523 \\
120  0.745 \\
130  1.173 \\
140  1.689 \\
150  2.304 \\
160  3.147 \\
170  4.395 \\
180  5.553 \\
190  7.214 \\
200  9.235 \\
};
%\addlegendentry{\normaliz};

\addplot [nmzstylealt]
  table[row sep=crcr]{
40    0.007 \\
50    0.014 \\
60    0.039 \\
70    0.080 \\
80    0.160 \\
90    0.300 \\
100   0.555 \\
110   0.931 \\
120   1.499 \\
130   2.342 \\
140   3.473 \\
150   5.120 \\
160   7.150 \\
170  10.116 \\
180  13.457 \\
190  18.045 \\
200  24.116 \\
};
%\addlegendentry{\normaliz *};

\end{axis}
\end{tikzpicture}%
}\hfill
  \subcaptionbox{$F_d(60)$\label{fig:knapsack:lattice:timings:plot:d}}{% This file was created by matlab2tikz v0.4.5 running on MATLAB 8.2.
% Copyright (c) 2008--2014, Nico Schlömer <nico.schloemer@gmail.com>
% All rights reserved.
% 
% The latest updates can be retrieved from
%   http://www.mathworks.com/matlabcentral/fileexchange/22022-matlab2tikz
% where you can also make suggestions and rate matlab2tikz.
% 
%
% defining custom colors
\definecolor{mycolor1}{rgb}{0.16863,0.50588,0.33725}%

\begin{tikzpicture}

% to keep aspect ratio
\setlength{\mywidth}{.4\linewidth} % you might want to change this
\setlength{\myheight}{.8\mywidth} % do not change this

\begin{axis}[%
width=\mywidth, %6.4cm,
height=\myheight, %5cm,
scale only axis,
xmin=4,
xmax=20,
ymode=log,
ymin=0.001,
ymax=10000,
yminorticks=true,
legend style={at={(0.2,0.02)},anchor=south west,draw=black,fill=white,legend cell align=left, font=\scriptsize, legend columns=4},
]

\addlegendimage{ftitstyle};
\addlegendentry{\ftit};
\addlegendimage{bboxstyle};
\addlegendentry{\bbox};
\addlegendimage{lattestyle};
\addlegendentry{\latte};
\addlegendimage{lattestylealt};
\addlegendentry{\latte{}$^\vee$};
\addlegendimage{nmzstyle};
\addlegendentry{\nmz};
\addlegendimage{nmzstylealt};
\addlegendentry{\nmz{}$^\vee$};
\addlegendimage{projstyle};
\addlegendentry{\proj};

\addplot [ftitstyle]
  table[row sep=crcr]{
4    7.722 \\
5   22.173 \\
6   33.795 \\
7   43.212 \\
8   44.388 \\
9   43.496 \\
10  43.577 \\
11  44.137 \\
12  44.183 \\
13  44.143 \\
14  44.275 \\
15  44.953 \\
16  44.967 \\
17  45.149 \\
18  45.115 \\
19  45.233 \\
20  45.396 \\
};
%\addlegendentry{\ftit};

\addplot [bboxstyle]
  table[row sep=crcr]{
4    0.236 \\ 
5    1.209 \\ 
6    3.861 \\ 
7    8.397 \\ 
8   18.088 \\
9   19.418 \\
10  20.616 \\
11  21.919 \\
12  23.121 \\
13  24.264 \\
14  25.290 \\
15  26.550 \\
16  27.808 \\
17  28.907 \\
18  30.089 \\
19  31.594 \\
20  32.797 \\
};	
%\addlegendentry{\bbox};

\addplot [lattestyle]
  table[row sep=crcr]{
4      0.007 \\   
5      0.015 \\   
6      0.030 \\   
7      0.077 \\   
8      0.188 \\   
9      0.507 \\   
10     1.315 \\   
11     3.521 \\   
12     9.238 \\   
13    24.864 \\  
14    66.012 \\  
15   177.424 \\ 
16   467.912 \\ 
17  1246.255 \\
18  3282.651 \\
};
%\addlegendentry{\latte};

\addplot [lattestylealt]
  table[row sep=crcr]{
4       0.009\\
5       0.027\\
6       0.113\\
7       0.504\\
8       2.120\\
9      11.207\\
10     44.764\\
11    255.896\\
12   2357.660\\
};
%\addlegendentry{\latte{}$^\vee$};

\addplot [projstyle]
  table[row sep=crcr]{
4    0.011\\
5    0.040\\
6    0.090\\
7    0.153\\
8    0.233\\
9    0.319\\
10   0.413\\
11   0.524\\
12   0.639\\
13   0.767\\
14   0.908\\
15   1.053\\
16   1.210\\
17   1.383\\
18   1.571\\
19   1.770\\
20   1.975\\
};
%\addlegendentry{\proj};

\addplot [nmzstyle]
  table[row sep=crcr]{
4    0.011 \\
5    0.193 \\
6    0.070 \\
7    0.148 \\
8    0.287 \\
9    0.365 \\
10   0.434 \\
11   0.502 \\
12   0.587 \\
13   0.689 \\
14  16.999 \\
15  18.982 \\
16  21.059 \\
17  22.979 \\
18  25.063 \\
19  27.296 \\
20  29.629 \\
};
%\addlegendentry{\nmz};

\addplot [nmzstylealt]
  table[row sep=crcr]{
4    0.019\\
5    0.036\\
6    0.047\\
7    0.054\\
8    0.054\\
9    0.060\\
10   0.064\\
11   0.066\\
12   0.069\\
13   0.073\\
14   0.074\\
15   0.077\\
16   0.082\\
17   0.084\\
18   0.084\\
19   0.090\\
20   0.092\\
};
%\addlegendentry{\nmz{}$^\vee$};

\end{axis}
\end{tikzpicture}%
}\\
  \caption{Timings (in seconds) for counting lattice points in Fibonacci knapsack polytopes $F_d(b)$ depending either on $d$ or on $b$, on a logarithmic scale. See Tables~\ref{tab:knapsack:lattice:rhs} and \ref{tab:knapsack:lattice:d} for exact timings.}
\label{fig:knapsack:lattice:timings:plot}
\end{figure}

We performed two series of experiments. In the first we fixed the dimension to $d=5$ and varied the right hand side of
the knapsack equation.  Geometrically speaking, this corresponds to shifting one of the facets of a simplex.  In the
second series we fixed the right hand side $b=60$ and varied the dimension. In both cases, we provided the full lists of
vertices and facets of the polytope as input to the algorithms, see~\cite{polymake_exp}.  All timings are listed in
Tables~\ref{tab:knapsack:lattice:rhs} and \ref{tab:knapsack:lattice:d} and plotted in
Figure~\ref{fig:knapsack:lattice:timings:plot}.  Note that in both cases the \projection algorithm performs pretty well.

We also determined what our memory restriction of 4 GB implies as an upper limit on the number of points depending on the dimension. This is corresponds roughly to 17~million points in $F_5(350)$ and 5~million points in $F_8(200)$. In both cases the \projection algorithm finishes in about a minute.

\subsubsection{Varying the right hand side}
First, note that when increasing the right hand side the number of lattice points $\ell$ in the polytope will increase
significantly.
The values of $\ell$ for our examples are also given in Table~\ref{tab:knapsack:lattice:rhs}.
Hence, the amount of work naturally increases for all algorithms but \latte.
More precisely, the running time for the Barvinok algorithm as implemented in \latte seems constant.
This is plausible, since changing the right hand side only gives a dilation of the same polytope, the amount of work to compute the generating functions at each vertex stays exactly the same.

This experiment, displayed in Figure~\ref{fig:knapsack:lattice:timings:plot:rhs}, shows the strength of Barvinok's algorithm, and this is why we suggest the Rule of Thumb~\ref{rot:lattice:latte}.

\subsubsection{Varying the dimension}
In this series of tests, we fix the right hand side, \textit{i.e.}, the size of the knapsack, and the coefficients, \textit{i.e.}, the object sizes, are the Fibonacci numbers.
At first sight the behavior of the implementations tested may seem erratic; see Figure~\ref{fig:knapsack:lattice:timings:plot:d}.
In fact, it correlates with various strategies to avoid the ``curse of dimension''.
If the dimension gets high enough the extra points will be bigger than the knapsack.
For our test this means that the number of lattice points, $\ell$, in the Fibonacci fractional knapsack polytopes becomes constant for $d \ge 8$; see Table~\ref{tab:knapsack:lattice:d}.
Moreover, all lattice points lie in the subspace corresponding to the first seven coordinates.

The lifting approach used in \ftit detects this and the running times are constant starting from $d=8$. The \bbox algorithm determines zero -- rounded down from some rational value less than one -- as the upper bound for each further coordinate and thus the only additional work comes from the increasing number of zero-entries in the lattice vectors.
The slight growth in the running time of the \projection algorithm mainly comes from the number of lifting steps that are needed to generate the lattice points.

Here, one can see that the dimension has a strong impact on the running time of \latte, because of the increasing number of vertices and the more complicated generating functions, both with all primal and default parameters. For high dimension, \ftit, \bbox, \normaliz, and \projection perform a lot better.
Since version~2.99.4, \normaliz employs a special strategy to reduce the size of the fundamental parallelepipeds, which is triggered above a certain threshold. This is the reason for the bumps at very low input parameters. The jump at dimension~14 is caused by an automatic coordinate upgrade that switches from 64 bit \texttt{long} arithmetic to arbitrary precision.
Further, the dual mode of \normaliz performs extremely well for this type of input.

\begin{ruleofthumb}
   For fractional knapsack polytopes try {\rm \projection}. For high dimensions also try {\rm \normaliz} in both default and dual mode.
\end{ruleofthumb}

\subsection{Random Box}
\label{sec:rbox}

We also considered a second class of polytopes that show a completely different behavior with respect to the different
algorithms. For parameters $d$ and $n$ we define a \emph{random box polytope} $R(d,n)$ as the convex hull of $n$ lattice
points chosen uniformly at random from the cube $C^d_5:=[0,5]^d$.

% REMOVED since 1st reviewer did not like it.
% temporarily creates unresolved references
%
% \begin{figure}[tb]
%   \setlength{\mywidth}{.27\linewidth} % you might want to change this
%   \subcaptionbox{\bbox (6,50)}{\input{img/histfit_rbox_bbox.tikz}}\qquad\qquad
%   \subcaptionbox{\projection (6,50)}{\input{img/histfit_rbox_proj.tikz}}\\
%   \subcaptionbox{\normaliz (6,50)}{\input{img/histfit_rbox_normaliz.tikz}}\qquad\qquad
%   \subcaptionbox{\latte (4,20)}{\input{img/histfit_rbox_latte.tikz}}\\
%   \caption{Histograms of timings of counting lattice points in 100 random polytopes, the boundaries of the $x$-axis correspond to the minimum and maximum values.}
%   \label{fig:rbox:stats:time}
% \end{figure}
%
%For each of the algorithms \bbox, \projection, \normaliz, and \latte we picked one combination of $(d,n)$ and applied
%the algorithm to $100$ randomly generated polytopes $R(d,n)$. The graphs in Figure~\ref{fig:rbox:stats:time} show the
%distributions of the running times of those algorithms for this pair of parameters.
%\todo[inline]{The histograms are missing. We dont have anything in the text about the statistical data, yet. Here is a small table:
%}
%\begin{itemize}
%\item \ftit does not work good on this. The average is: $10796.377$. Standard Deviation is: $65314.382$
%\item \bbox Average: $31.165$. Standard Deviation is: $4.34$
%\item \latte Average: $1.601$. Standard Deviation is: $0.560$
%\item \normaliz Average: $0.136$. Standard Deviation is: $0.013$
%\item \projection Average: $6.2$. Standard Deviation is: $0.853$
%\end{itemize}
%\todo[inline]{end table}

\begin{figure}[p]
  \centering
  \subcaptionbox{\bbox}{% This file was created by matlab2tikz v0.4.7 (commit 1a3c9d1cf3aa49ad44d17ce5def556346dfbdc20) running on MATLAB 8.3.
% Copyright (c) 2008--2014, Nico Schlömer <nico.schloemer@gmail.com>
% All rights reserved.
% Minimal pgfplots version: 1.3
% 
% The latest updates can be retrieved from
%   http://www.mathworks.com/matlabcentral/fileexchange/22022-matlab2tikz
% where you can also make suggestions and rate matlab2tikz.
% 
\begin{tikzpicture}

% to keep aspect ratio
\setlength{\mywidth}{.33\linewidth} % you might want to change this
\setlength{\myheight}{.85\mywidth} % do not change this

\begin{axis}[%
width=\mywidth,%5.15805555555556in,
height=\myheight,%4.34666666666667in,
view={-37.5}{30},
scale only axis,
xmin=20,
xmax=70,
xmajorgrids,
ymin=4,
ymax=9,
ymajorgrids,
zmin=0,
zmax=4000,
zmajorgrids,
axis x line*=bottom,
axis y line*=left,
axis z line*=left,
font=\footnotesize,
ticklabel style={/pgf/number format/.cd,1000 sep={\,},}
]
\addplot3[%
mesh,
very thick,
shader=flat,
mesh/rows=6]
table[row sep=crcr,header=false] {
20	4	   0.082\\
20	5	   0.756\\
20	6	   4.298\\
20	7	  20.466\\
20	8	  89.321\\
20	9	 473.808\\
30	4	   0.130\\
30	5	   1.439\\
30	6	  11.853\\
30	7	  52.392\\
30	8	 270.244\\
30	9	1005.608\\
40	4	   0.154\\
40	5	   2.548\\
40	6	  18.827\\
40	7	 144.850\\
40	8	 679.447\\
40	9	3511.637\\
50	4	   0.190\\
50	5	   3.116\\
50	6	  34.396\\
50	7	 236.414\\
50	8	1693.210\\
50	9	3600    \\
60	4	   0.207\\
60	5	   3.766\\
60	6	  47.261\\
60	7	 400.435\\
60	8	2929.315\\
60	9	3600    \\
70	4	   0.213\\
70	5	   4.845\\
70	6	  63.581\\
70	7	 581.584\\
70	8	3600    \\
70	9	3600    \\
};
\end{axis}
\end{tikzpicture}%
}
  \hfill
  \subcaptionbox{\projection}{% This file was created by matlab2tikz v0.4.7 (commit 1a3c9d1cf3aa49ad44d17ce5def556346dfbdc20) running on MATLAB 8.3.
% Copyright (c) 2008--2014, Nico Schlömer <nico.schloemer@gmail.com>
% All rights reserved.
% Minimal pgfplots version: 1.3
% 
% The latest updates can be retrieved from
%   http://www.mathworks.com/matlabcentral/fileexchange/22022-matlab2tikz
% where you can also make suggestions and rate matlab2tikz.
% 
\begin{tikzpicture}

% to keep aspect ratio
\setlength{\mywidth}{.33\linewidth} % you might want to change this
\setlength{\myheight}{.85\mywidth} % do not change this

\begin{axis}[%
width=\mywidth,%5.15805555555556in,
height=\myheight,%4.34666666666667in,
view={-37.5}{30},
scale only axis,
xmin=20,
xmax=70,
xmajorgrids,
ymin=4,
ymax=9,
ymajorgrids,
zmin=0,
zmax=4000,
zmajorgrids,
axis x line*=bottom,
axis y line*=left,
axis z line*=left,
font=\footnotesize,
ticklabel style={/pgf/number format/.cd,1000 sep={\,},}
]
\addplot3[%
mesh,
very thick,
shader=flat,
mesh/rows=6]
table[row sep=crcr,header=false] {
20	4	   0.013 \\
20	5	   0.100 \\
20	6	   0.450 \\
20	7	   1.369 \\
20	8	   3.414 \\
20	9	   7.490 \\
30	4	   0.023 \\
30	5	   0.232 \\
30	6	   1.806 \\
30	7	   9.502 \\
30	8	  41.705 \\
30	9	 151.081 \\
40	4	   0.029 \\
40	5	   0.462 \\
40	6	   3.945 \\
40	7	  28.062 \\
40	8	 158.855 \\
40	9	 219.040 \\
50	4	   0.035 \\
50	5	   0.621 \\
50	6	   6.306 \\
50	7	  57.688 \\
50	8	 138.180 \\
50	9	 248.157 \\
60	4	   0.042 \\
60	5	   0.762 \\
60	6	  10.476 \\
60	7	 102.579 \\
60	8	 148.227 \\
60	9	 769.460 \\
70	4	   0.040 \\
70	5	   1.032 \\
70	6	  14.247 \\
70	7	 149.536 \\
70	8	 223.361 \\
70	9	1017.670 \\
};
\end{axis}
\end{tikzpicture}%
} \\
  \subcaptionbox{\label{fig:rbox:lattice:latte}\latte}{% This file was created by matlab2tikz v0.4.7 (commit 1a3c9d1cf3aa49ad44d17ce5def556346dfbdc20) running on MATLAB 8.3.
% Copyright (c) 2008--2014, Nico Schlömer <nico.schloemer@gmail.com>
% All rights reserved.
% Minimal pgfplots version: 1.3
% 
% The latest updates can be retrieved from
%   http://www.mathworks.com/matlabcentral/fileexchange/22022-matlab2tikz
% where you can also make suggestions and rate matlab2tikz.
% 
\begin{tikzpicture}

% to keep aspect ratio
\setlength{\mywidth}{.33\linewidth} % you might want to change this
\setlength{\myheight}{.85\mywidth} % do not change this

\begin{axis}[%
width=\mywidth,%5.15805555555556in,
height=\myheight,%4.34666666666667in,
view={-37.5}{30},
scale only axis,
xmin=20,
xmax=70,
xmajorgrids,
ymin=4,
ymax=9,
ymajorgrids,
zmin=0,
zmax=4000,
zmajorgrids,
axis x line*=bottom,
axis y line*=left,
axis z line*=left,
font=\footnotesize,
ticklabel style={/pgf/number format/.cd,1000 sep={\,},}
]
\addplot3[%
mesh,
very thick,
shader=flat,
mesh/rows=6]
table[row sep=crcr,header=false] {
20	4	  1.798 \\
20	5	438.684 \\
20	6	3600    \\
20	7	3600    \\
20	8	3600    \\
20	9	3600    \\
30	4	  1.639 \\
30	5	612.918 \\
30	6	3600    \\
30	7	3600    \\
30	8	3600    \\
30	9	3600    \\
40	4	  1.419 \\
40	5	534.548 \\
40	6	3600    \\
40	7	3600    \\
40	8	3600    \\
40	9	3600    \\
50	4	  1.357 \\
50	5	588.029 \\
50	6	3600    \\
50	7	3600    \\
50	8	3600    \\
50	9	3600    \\
60	4	  1.405 \\
60	5	507.007 \\
60	6	3600    \\
60	7	3600    \\
60	8	3600    \\
60	9	3600    \\
70	4	  1.172 \\
70	5	404.041 \\
70	6	3600    \\
70	7	3600    \\
70	8	3600    \\
70	9	3600    \\
};
\end{axis}
\end{tikzpicture}%
}
  \hfill
  \subcaptionbox{\label{fig:rbox:lattice:lattevee}\latte{}$^\vee$}{% This file was created by matlab2tikz v0.4.7 (commit 1a3c9d1cf3aa49ad44d17ce5def556346dfbdc20) running on MATLAB 8.3.
% Copyright (c) 2008--2014, Nico Schlömer <nico.schloemer@gmail.com>
% All rights reserved.
% Minimal pgfplots version: 1.3
% 
% The latest updates can be retrieved from
%   http://www.mathworks.com/matlabcentral/fileexchange/22022-matlab2tikz
% where you can also make suggestions and rate matlab2tikz.
% 
\begin{tikzpicture}

% to keep aspect ratio
\setlength{\mywidth}{.33\linewidth} % you might want to change this
\setlength{\myheight}{.85\mywidth} % do not change this

\begin{axis}[%
width=\mywidth,%5.15805555555556in,
height=\myheight,%4.34666666666667in,
view={-37.5}{30},
scale only axis,
xmin=20,
xmax=70,
xmajorgrids,
ymin=4,
ymax=9,
ymajorgrids,
zmin=0,
zmax=4000,
zmajorgrids,
axis x line*=bottom,
axis y line*=left,
axis z line*=left,
font=\footnotesize,
ticklabel style={/pgf/number format/.cd,1000 sep={\,},}
]
\addplot3[%
mesh,
very thick,
shader=flat,
mesh/rows=6]
table[row sep=crcr,header=false] {
20	4	0.252    \\
20	5	3.111    \\
20	6	49.751   \\
20	7	801.512  \\
20	8	3600     \\
20	9	3600     \\
30	4	0.375    \\
30	5	4.651    \\
30	6	122.619  \\
30	7	3016.324 \\
30	8	3600     \\
30	9	3600     \\
40	4	0.49     \\
40	5	6.69     \\
40	6	154.804  \\
40	7	3600     \\
40	8	3600     \\
40	9	3600     \\
50	4	0.512    \\
50	5	8.114    \\
50	6	235.039  \\
50	7	3600     \\
50	8	3600     \\
50	9	3600     \\
60	4	0.569    \\
60	5	9.516    \\
60	6	281.204  \\
60	7	3600     \\
60	8	3600     \\
60	9	3600     \\
70	4	0.672    \\
70	5	11.017   \\
70	6	338.06   \\
70	7	3600     \\
70	8	3600     \\
70	9	3600     \\
};
\end{axis}
\end{tikzpicture}%
}\\
  \subcaptionbox{\normaliz}{% This file was created by matlab2tikz v0.4.7 (commit 1a3c9d1cf3aa49ad44d17ce5def556346dfbdc20) running on MATLAB 8.3.
% Copyright (c) 2008--2014, Nico Schlömer <nico.schloemer@gmail.com>
% All rights reserved.
% Minimal pgfplots version: 1.3
% 
% The latest updates can be retrieved from
%   http://www.mathworks.com/matlabcentral/fileexchange/22022-matlab2tikz
% where you can also make suggestions and rate matlab2tikz.
% 
\begin{tikzpicture}

% to keep aspect ratio
\setlength{\mywidth}{.33\linewidth} % you might want to change this
\setlength{\myheight}{.85\mywidth} % do not change this

\begin{axis}[%
width=\mywidth,%5.15805555555556in,
height=\myheight,%4.34666666666667in,
view={-37.5}{30},
scale only axis,
xmin=20,
xmax=70,
xmajorgrids,
ymin=4,
ymax=9,
ymajorgrids,
zmin=0,
zmax=4000,
zmajorgrids,
axis x line*=bottom,
axis y line*=left,
axis z line*=left,
font=\footnotesize,
ticklabel style={/pgf/number format/.cd,1000 sep={\,},}
]
\addplot3[%
mesh,
very thick,
shader=flat,
mesh/rows=6]
table[row sep=crcr,header=false] {
20	4	  0.002 \\
20	5	  0.005 \\
20	6	  0.017 \\
20	7	  0.098 \\
20	8	  1.068 \\
20	9	  5.256 \\
30	4	  0.003 \\
30	5	  0.007 \\
30	6	  0.040 \\
30	7	  0.426 \\
30	8	  3.214 \\
30	9	 50.533 \\
40	4	  0.004 \\
40	5	  0.013 \\
40	6	  0.085 \\
40	7	  0.782 \\
40	8	  8.565 \\
40	9	129.305 \\
50	4	  0.006 \\
50	5	  0.020 \\
50	6	  0.129 \\
50	7	  1.307 \\
50	8	 16.734 \\
50	9	234.613 \\
60	4	  0.006 \\
60	5	  0.022 \\
60	6	  0.175 \\
60	7	  1.941 \\
60	8	 24.145 \\
60	9	474.163 \\
70	4	  0.005 \\
70	5	  0.025 \\
70	6	  0.237 \\
70	7	  2.392 \\
70	8	 34.041 \\
70	9	626.180 \\
};
\end{axis}
\end{tikzpicture}%
}
%  \hfill
%  \subcaptionbox{\normaliz{}$^\vee$}{\input{img/rbox_lattice_plot_normaliz_flags.tikz}}

  \caption{Running times (in seconds) for counting lattice points in random polytopes depending on the dimension $d\in[4,9]$ and the number of points $m\in[20,70]$. See Table~\ref{tab:randombox:integertimings:} for exact timings.}
  \label{fig:rbox:lattice:timings}
\end{figure}
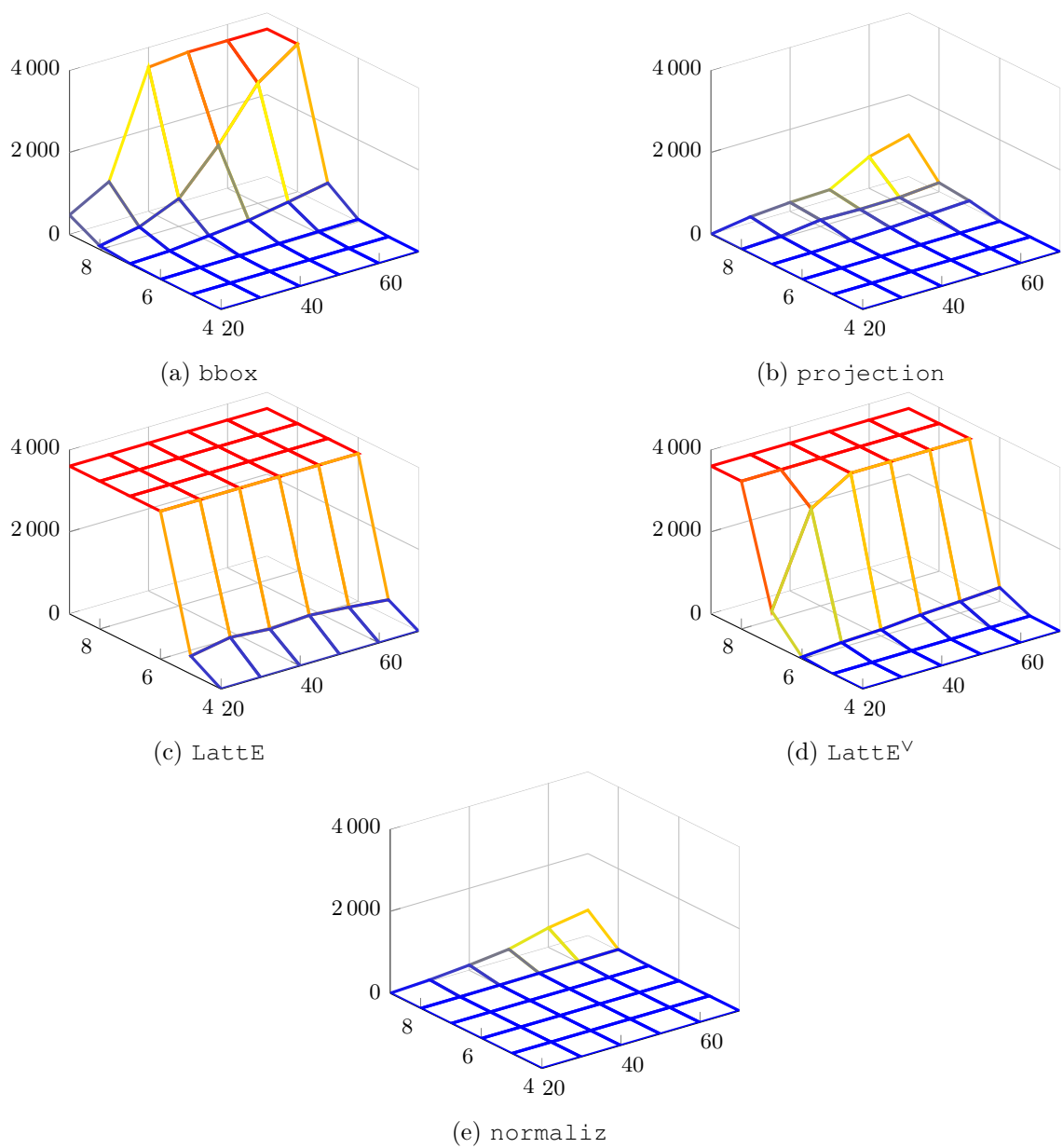

For each of the algorithms \bbox, \latte, \normaliz, and \projection we did the following experiment: for each dimension
$d$ between $4$ and $9$, and for each $n$ in the set $\{20,30,40,50,60,70\}$ we generated $10$ random polytopes $R(d,n)$
and tried to count the lattice points in the convex hull with each of the algorithms. Again, we provided the full lists
of vertices and facets of the polytope as input to the algorithms, see~\cite{polymake_exp}. The average running times
are plotted in Figure~\ref{fig:rbox:lattice:timings} and are listed in Table~\ref{tab:randombox:integertimings:}.  We
did not include \ftit and \normaliz in dual mode in these experiments; for details see Remark~\ref{rem:rbox:ftitnmz}.

In this experiment only \normaliz and \projection finished all tests within the given
time and memory restrictions.
%Note that the dual algorithm of \normaliz does not behave well with this kind of input,
%which fits the authors recommendation.
\bbox failed on some of the $9$-dimensional examples, while \latte basically only worked for
dimension~$4$ and~$5$. With the \emph{all primal} mode, \latte was able to compute some of the $7$-dimensional examples.

Observe that the number of lattice points inside $C^d_5$ grows exponential in $d$ (in fact, there are $6^d$ lattice
points). Hence, the fraction of the lattice points chosen decreases with the dimension.
Yet, we expect that coordinate projections of the polytope $R(d,n)$ use the
full range $[0,5]$ for its coordinates, which forces \bbox to enumerate many lattice points that do not
contribute to the final result. This clarifies why this algorithm runs into problems on this class of polytopes for higher dimensions.

We see from the experiments that \latte basically only works for the low dimensional examples.
This particularly bad behavior can be explained by the specific version of Barvinok's algorithm employed.

We have already outlined the basic idea of this algorithm above.  In its original form it starts with some initial
triangulation of the homogenization of the polytope into simplicial cones. The index of such a simplicial cone is the
determinant of its generators. The algorithm now successively lowers the maximum index of a cone in the triangulation by
replacing cones with a signed decomposition into cones of lower index. The number of steps used by the algorithm heavily
depends on the initial maximal index of a cone. Once we have a signed decomposition we can add up the rational
generating functions for each cone. A drawback of this approach is that we count lattice points in intersections of
cones multiple times, once for each cone they are in. \latte implements both options discussed above to avoid this
problem. We can choose to work with the dual cone, or we can use \emph{irrational decomposition}.

By default, \latte chooses the first of the two variants and does a triangulation in dual space. For our experiments we
can now observe that the index of an initial triangulation of the dual of a random polytope $R(d,n)$ is quite
large. This explains that \latte won't finish on these examples within our given time limit.

We observe, however, that the index of an initial triangulation in primal space is much lower. So the performance of
\latte significantly improves on these examples if we ask \latte to use \emph{irrational decomposition} and do a signed
triangulation in primal space; compare Figures~\ref{fig:rbox:lattice:latte} and~\ref{fig:rbox:lattice:lattevee}.

Since we are dealing with random polytopes we did a similar statistical analysis as at the end of Section~\ref{sec:voronoi}, to explain why it suffices to list the average over only ten runs per parameter set.
The natural random variables are the running times for the various algorithms.
For \bbox, \latte, \normaliz, and \projection, they are strongly concentrated around their respective mean.
We omit the details.

\begin{figure}[ht]
  \centering
  \setlength{\mywidth}{.27\linewidth} % you might want to change this
  %\subcaptionbox{\ftit (4,20)}{\input{img/histfit_rbox_4ti2.tikz}}\qquad\qquad
  \subcaptionbox{\ftit (4,20)}{% This file was created by matlab2tikz v0.4.5 running on MATLAB 8.2.
% Copyright (c) 2008--2014, Nico Schlömer <nico.schloemer@gmail.com>
% All rights reserved.
% 
% The latest updates can be retrieved from
%   http://www.mathworks.com/matlabcentral/fileexchange/22022-matlab2tikz
% where you can also make suggestions and rate matlab2tikz.
% 
%
% defining custom colors
\definecolor{mycolor1}{rgb}{0.16863,0.50588,0.33725}%
\begin{tikzpicture}

% to keep aspect ratio
%\setlength{\mywidth}{.3\linewidth} % you might want to change this
\setlength{\myheight}{.8\mywidth} % do not change this

\begin{axis}[%
width=\mywidth,
height=\myheight,
scale only axis,
xmin=1,
xmax=90,%100,
ymode=log,
ymin=0.1,
ymax=1000000,
yminorticks=true,
legend style={empty legend}%{at={(0.74,0.02)},anchor=south west,draw=black,fill=white,legend cell align=left}
]
\addplot [color=red,draw=none,solid,mark=*,mark options={solid,fill=red}]
  table[row sep=crcr]{
1     9.4999999999   \\
2     38.03          \\
3     704.33         \\
4     28.78          \\
5     51.42          \\
6     47.590000000   \\
7     37.99          \\
8     5.4799999999   \\
9     16.749999999   \\
10    9.8100000000   \\
11    1.3200000000   \\
12    2011.35        \\
13    1170.15        \\
14    45.789999999   \\
15    458977.51      \\
16    3.1300000000   \\
17    19.809999999   \\
18    19.960000000   \\
19    2.1899999999   \\
20    154.13000000   \\
21    6.3999999999   \\
22    0.7200000000   \\
23    1.9199999999   \\
24    633.96999999   \\
25    6.7300000000   \\
26    2.0300000000   \\
27    23.869999999   \\
28    634.97999999   \\
29    4.8300000000   \\
30    693.93999999   \\
31    311.62999999   \\
32    1.5600000000   \\
33    23.220000000   \\
34    1776.9299999   \\
35    29.990000000   \\
36    2.7999999999   \\
37    16.969999999   \\
38    3.5199999999   \\
39    7.2499999999   \\
40    5.7900000000   \\
41    1.2600000000   \\
42    11541.33       \\
43    10.280000000   \\
44    2.9299999999   \\
45    1.6500000000   \\
46    1.8000000000   \\
47    5.5899999999   \\
48    43.550000000   \\
49    22.069999999   \\
50    425.28000000   \\
51    1500.84        \\
52    3239.1299999   \\
53    14.860000000   \\
54    8.5200000000   \\
55    4.1699999999   \\
56    45.62          \\
57    57.290000000   \\
58    5.9099999999   \\
59    505.11000000   \\
60    2.6299999999   \\
61    1.6499999999   \\
62    123.18000000   \\
63    171.73         \\
64    34.350000000   \\
65    15.319999999   \\
66    5.1600000000   \\
67    1.1599999999   \\
68    92.349999999   \\
69    4.2900000000   \\
70    4.0799999999   \\
71    1.2499999999   \\
72    6.6700000000   \\
73    18.980000000   \\
74    119.65         \\
75    269.44999999   \\
76    52726.720000   \\
77    11.939999999   \\
78    1.7699999999   \\
79    10.740000000   \\
80    2.3300000000   \\
81    4.4299999998   \\
82    1.2600000000   \\
83    601.52999999   \\
84    1114.3900000   \\
85    3.1900000000   \\
86    118.67000000   \\
87    409631.05      \\
};
%\addlegendentry{\ftit};

\end{axis}
\end{tikzpicture}%
}\qquad\qquad
  \subcaptionbox{\normaliz{}$^\vee$ (4,20)}{% This file was created by matlab2tikz v0.4.5 running on MATLAB 8.2.
% Copyright (c) 2008--2014, Nico Schlömer <nico.schloemer@gmail.com>
% All rights reserved.
% 
% The latest updates can be retrieved from
%   http://www.mathworks.com/matlabcentral/fileexchange/22022-matlab2tikz
% where you can also make suggestions and rate matlab2tikz.
% 
%
% defining custom colors
\definecolor{mycolor1}{rgb}{0.16863,0.50588,0.33725}%
\begin{tikzpicture}

% to keep aspect ratio
%\setlength{\mywidth}{.3\linewidth} % you might want to change this
\setlength{\myheight}{.8\mywidth} % do not change this

\begin{axis}[%
width=\mywidth,
height=\myheight,
scale only axis,
xmin=1,
xmax=90,
ymode=log,
ymin=0.005,
ymax=1000,
yminorticks=true,
legend style={empty legend}%{at={(0.74,0.02)},anchor=south west,draw=black,fill=white,legend cell align=left}
]
\addplot [color=red,draw=none,solid,mark=*,mark options={solid,fill=red}]
  table[row sep=crcr]{
1     0.0600000000000023    \\
2     0.570000000000002     \\
3     0.379999999999995     \\
4     0.0199999999999818    \\
5     0.189999999999998     \\
6     0.0499999999999829    \\
7     0.0400000000000205    \\
8     0.0800000000000125    \\
9     0.0400000000000205    \\
10    0.0300000000000011    \\
11    0.00999999999999091   \\
12    0.0199999999999818    \\
13    0.120000000000005     \\
14    0.71999999999998      \\
15    0.0600000000000023    \\
16    0.0200000000000102    \\
17    0.110000000000014     \\
18    0.0300000000000011    \\
19    0.0199999999999818    \\
20    0.00999999999999091   \\
21    0.0300000000000011    \\
22    0.0200000000000102    \\
23    0.0199999999999818    \\
24    55.46                 \\
25    0.109999999999985     \\
26    0.0400000000000205    \\
27    0.0500000000000114    \\
28    0.0300000000000011    \\
29    0.0699999999999932    \\
30    0.0400000000000205    \\
31    7.86000000000001      \\
32    0.039999999999992     \\
33    0.0199999999999818    \\
34    0.449999999999989     \\
35    0.0700000000000216    \\
36    0.0300000000000011    \\
37    0.0199999999999818    \\
38    0.0700000000000216    \\
39    0.00999999999999091   \\
40    0.039999999999992     \\
41    0.0199999999999818    \\
42    0.109999999999985     \\
43    0.0400000000000205    \\
44    2.44999999999999      \\
45    0.0200000000000102    \\
46    0.0300000000000011    \\
47    0.0300000000000011    \\
48    0.0200000000000102    \\
49    0.110000000000014     \\
50    0.0200000000000102    \\
51    2.96000000000001      \\
52    0.0400000000000205    \\
53    0.0299999999999727    \\
54    0.109999999999957     \\
55    0.0699999999999932    \\
56    0.0199999999999818    \\
57    0.0299999999999727    \\
58    0.0199999999999818    \\
59    0.20999999999998      \\
60    0.00999999999999091   \\
61    0.0100000000000477    \\
62    0.0400000000000205    \\
63    0.0500000000000114    \\
64    0.0799999999999841    \\
65    0.0299999999999727    \\
66    0.0299999999999727    \\
67    0.0299999999999727    \\
68    0.0299999999999727    \\
69    0.0600000000000023    \\
70    0.0699999999999932    \\
71    0.0299999999999727    \\
72    52.85                 \\
73    0.0400000000000205    \\
74    0.539999999999964     \\
75    0.0500000000000114    \\
76    0.740000000000018     \\
77    0.07000000000005      \\
78    0.00999999999999091   \\
79    0.0400000000000205    \\
80    0.110000000000014     \\
81    0.0200000000000387    \\
82    0.0100000000000477    \\
83    0.0299999999999727    \\
84    0.0600000000000023    \\
85    0.0199999999999818    \\
86    0.0199999999999818    \\
87    0.0300000000000296    \\
};

% \addlegendentry{\nmz{}$^\vee$};

\end{axis}
\end{tikzpicture}%
}
%  \caption{Histogram of timings of counting lattice points in random polytopes using
%    \ftit, together with a plot of the individual running times. The
%    figures show that \ftit is extremely sensitive to the actual input.}
  \caption{Plot of the running times of counting lattice points in 87 random polytopes using \ftit and \normaliz in dual mode. The figures show that this kind of algorithm is very sensitive to the actual input.}
  \label{fig:rbox:stats:time:4ti2nmz}
\end{figure}
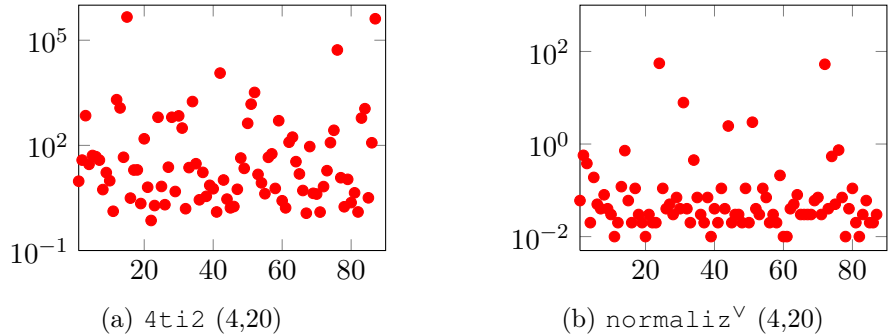

\begin{remark}\label{rem:rbox:ftitnmz}
We omit the timings for \ftit and \normaliz in dual mode because for this experiment the statistical analysis found that they behave erratically, \normaliz to a lesser extent than \ftit; see Figure~\ref{fig:rbox:stats:time:4ti2nmz}. We do not have an explanation for this phenomenon.
\end{remark}

\subsection{Matching Polytopes}
\label{sec:matching}

Consider a finite simple undirected graph $G=(V,E)$.
For a vertex $v\in V$ we let $\delta(v)\subseteq E$ be the set of edges incident with $v$.
A vector $x\in \RR^E$ which satisfies
  \[
  \begin{array}{r@{\quad}l}
    \sum_{e\in\delta(v)} x_e \le 1 & \text{for all } v\in V \quad \text{and} \\
    x_e\ge 0 & \text{for all } e\in E
  \end{array}
  \]
is called a \emph{fractional matching}. This gives rise to the \emph{fractional matching polytope}
  \[
    M(G) \ := \ \conv(x\in \RR^E \mid x \text{ is a fractional matching of } G)
  \]
of the graph $G$.
The dimension $d$ of $M(G)$ is easy to determine, it is the same as the number of edges $d=|E|$.
Further the number $m$ of constraints on $M(G)$ is equal to $|V| + |E|$.

A \emph{matching} of $G$ is a subset of the edge set $E$ which covers each vertex at most once.
By construction the matchings are in bijection with the integer points in $M(G)$, all of which have $0/1$-coordinates.
Finding matchings, e.g., of maximum cardinality, has a variety of applications.
Thus, it is generally desirable to solve integer linear programs over $M(G)$.
This is our motivation to look into counting and enumerating matchings, i.e., the lattice points in fractional matching polytopes.

\begin{figure}[tb]
  % This file was created by matlab2tikz v0.4.5 running on MATLAB 8.2.
% Copyright (c) 2008--2014, Nico Schlömer <nico.schloemer@gmail.com>
% All rights reserved.
% 
% The latest updates can be retrieved from
%   http://www.mathworks.com/matlabcentral/fileexchange/22022-matlab2tikz
% where you can also make suggestions and rate matlab2tikz.
% 
%
% defining custom colors
\definecolor{mycolor1}{rgb}{0.16863,0.50588,0.33725}%
\begin{tikzpicture}

% to keep aspect ratio
\setlength{\mywidth}{.4\linewidth} % you might want to change this
\setlength{\myheight}{.8\mywidth} % do not change this

\begin{axis}[%
width=\mywidth, %6.4cm,
height=\myheight, %5cm,
scale only axis,
xmin=4,
xmax=16,
ymode=log,
ymin=0.001,
ymax=10000,
yminorticks=true,
legend style={at={(0.55,1.02)},anchor=north west,draw=black,fill=white,legend cell align=left, font=\scriptsize, legend columns=2}
]

\addlegendimage{ftitstyle};
\addlegendentry{\ftit};
\addlegendimage{azovestyle};
\addlegendentry{\azove};
\addlegendimage{bboxstyle};
\addlegendentry{\bbox};
\addlegendimage{lattestyle};
\addlegendentry{\latte};
\addlegendimage{lattestylealt};
\addlegendentry{\latte{}$^\vee$};
\addlegendimage{nmzstyle};
\addlegendentry{\nmz};
\addlegendimage{nmzstylealt};
\addlegendentry{\nmz{}$^\vee$};
\addlegendimage{projstyle};
\addlegendentry{\proj};

\addplot [ftitstyle]
  table[row sep=crcr]{
4      0.004 \\
5      0.004 \\
6      0.008 \\
7      0.184 \\
8     12.058 \\
9   1557.383 \\
};
%\addlegendentry{\ftit};

\addplot [bboxstyle]
  table[row sep=crcr]{
4      0.005 \\ 
5      0.051 \\ 
6      2.430 \\ 
7    281.244 \\ 
8   6138.776 \\
};	
%\addlegendentry{\bbox};

\addplot [lattestyle]
  table[row sep=crcr]{
4        0.028\\   
5        0.450\\   
6       13.195\\   
7      531.014\\   
};
%\addlegendentry{\latte};

\addplot [lattestylealt]
  table[row sep=crcr]{
4        0.07\\   
5       16.865\\   
};
%\addlegendentry{\latte{}$^\vee$};

\addplot [projstyle]
  table[row sep=crcr]{
4     0.002\\
5     0.009\\
6     0.059\\
7     0.507\\
8     4.306\\
9    60.759\\
};
%\addlegendentry{\projection};

\addplot [nmzstyle]
  table[row sep=crcr]{
4      0.001 \\
5      0.006 \\
6      0.093 \\
7   1782.251 \\
};
%\addlegendentry{\normaliz};

\addplot [nmzstylealt]
  table[row sep=crcr]{
4    0.001 \\
5    0.003 \\
6    0.008 \\
7    0.072 \\
8    2.163 \\
9   95.743 \\
10 558.374 \\
};
%\addlegendentry{\normaliz{}^$\vee$};

\addplot [azovestyle]
  table[row sep=crcr]{
4    0.001 \\
5    0.001 \\
6    0.003 \\
7    0.004 \\
8    0.004 \\
9    0.004 \\
10   0.009 \\
11   0.017 \\
12   0.062 \\
13   0.236 \\
14   1.046 \\
15   4.696 \\
16  23.224 \\
};
%\addlegendentry{\azove};

\end{axis}
\end{tikzpicture}%
  \centering
  \caption{Timings (in seconds) of counting lattice points in the matching polytope of $K_n$, on a logarithmic scale; see Table~\ref{tab:matching:lattice} for exact timings.}
  \label{fig:matching:lattice}
\end{figure}
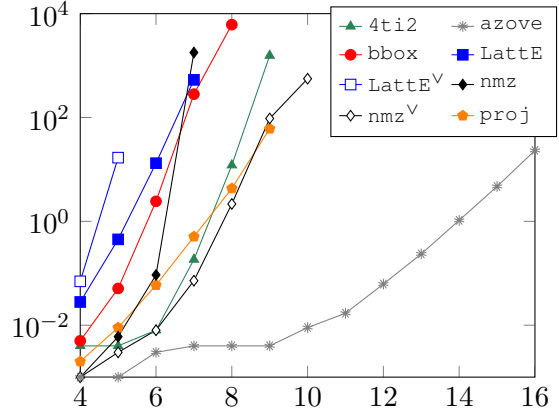

Our test graph is the complete graph $K_n$ with $n$ nodes.  In this case $d = \binom{n}{2}$ is the dimension of
$M(K_n)$, and the number of inequalities equals $m = \binom{n}{2} + n$.  In addition to the five implementations used
before, we also tested \azove.  This is a special implementation for enumerating $0/1$-vectors.  Each code was tried ten
times on each input with the same parameter. In each experiment, the polytopes were given by their list of facets and,
for the \projection algorithm, also the vertices, see~\cite{polymake_exp}.

The average running times are displayed in Figure~\ref{fig:matching:lattice} and Table~\ref{tab:matching:lattice}.

Note that \projection could go slightly further but this algorithm does not only require the inequalities of the
polytope but also the vertices.  As $M(K_9)$ already has as many as $79\,892$ vertices it was not possible in our setup
to compute the vertices of $M(K_{10})$: we expect around $700\,000$ vertices which exceed our memory constraint of 4~GB.

Since we start with an outer description, the dual algorithm of \normaliz again performs significantly better.

It is no surprise that \azove outperforms all other methods, as it was specifically designed to solve this kind of problem.
No other code can obtain the result in less than one hour for $n\ge 10$.
For comparison, it takes $2916.652$ seconds for \azove to finish the task for $M(K_{19})$.

\begin{ruleofthumb}
  Use {\rm \azove} to enumerate the $0/1$-points in a polytope.
\end{ruleofthumb}

Notice, however, that \azove might need to be patched; see Remark~\ref{rem:azove}.

\subsection{Additional Comments}
\label{sec:lattice:additional_comments}

The above examples show that the performance of the various algorithms does depend on properties of the polytopes that
are difficult to check a priori. The following example also shows that a high initial index need not necessarily lead to problems for \latte. Let
\begin{align*}
  P\ :=\ \conv(0,e_1,e_2,e_1+e_2+a\cdot e_3, e_1+e_2+b\cdot e_4, e_1+e_2+c\cdot e_5)
\end{align*}
for large pairwise co-prime $a,b,c$. This polytope is a simplex whose only lattice points are the vertices. Its index is
$abc$. Still \latte finishes fast, \projection at least finishes, while \bbox, \ftit, and \normaliz either do not
complete the computation within reasonable time or run into numerical problems.

This example also demonstrates that the problem of counting or enumerating lattice points in polytopes has a distinctive
number theoretic flavor. Small changes in the input which, from a linear programming point of view, have only little
impact on the problem may behave completely different with respect to the counting problem, as small changes may greatly
influence the number of \emph{integer solutions} of an equation. The algorithms used in \ftit, \latte, and \normaliz all
implicitly depend on such number theoretic influences, which explains some of the results we observed. In contrast,
the \projection method implemented in \polymake uses a plain convex geometric approach. This suggests why the latter
approach is usually not the best on any given input but often behaves quite well.
\begin{ruleofthumb}
  If you want to enumerate lattice points in a general polytope and do not have any further information then try {\rm \projection}.
\end{ruleofthumb}

\section{Summary}
\label{sec:summary}

\noindent
For our discussion on polyhedral and optimization algorithms and their implementations we have focused on two key
problems. The \emph{representation conversion problem} is a central task in classical polyhedral geometry, while the
problem of \emph{counting} or \emph{enumerating} lattice points is a main task in integer or combinatorial optimization,
number theory and algebra. Commonly used algorithms and implementations for these two problems are conveniently
available through \polymake. For each of the problems \polymake itself adds additional algorithms, the \emph{beneath-and-beyond
method} for the conversion problem and both the \emph{projection} as well as the \emph{bounding box method} for the enumeration of lattice points.

We have explained with various examples that in particular our implementation for the projection method behaves
surprisingly well on a wide range of problems, while many of the other algorithms show a much stronger dependence on the
type of the input.  For the convex hull algorithms the picture is less clear as each algorithm and each implementation
has its particular class of problems where it is most successful.

In our Rule of Thumb~\ref{rule:double-description} we advertised the double description method.
Many of our tests seem to suggest that its implementation in \ppl is superior to the one in \cdd.
It should be noted, however, that there are instances known where that superiority is reversed.
For instance, Avis and Jordan \cite{AvisJordan:1510.02545} report on cases where \ppl is a hundred times slower than \cdd.

For lack of maintenance the software packages \porta and \azove are essentially outdated, but they can still sometimes be useful.
However, any naive usage may result in crashes.

Our examples were chosen with the desire to look into a small number of plausible scenarios.
The main message is that very often it is difficult to predict which convex hull algorithm performs best on a specific input.
Nonetheless, we hope that our ten ``rules of thumb'' can serve as a guideline.

\appendix
\setcounter{table}{0}
\renewcommand{\thetable}{\thesection.\arabic{table}}
\section{Experimental Setup}
\label{sec:setup}

\noindent
Everything was calculated on identical Linux machines with the memory limit set to 4~GB (via \texttt{ulimit}).
All tests in one section were done on the same machine.
Any test exceeding this bound is marked as \textoom in the respective tables.
All timings were measured in CPU seconds, except for tests on non-symmetric cut polytopes where we used wallclock time to show the performance of the multithreaded version of \normaliz.
Those entries marked with a ${}^*$ ran only one iteration.
The hardware for all tests was:
\begin{quote}
  CPU: AMD Phenom(tm) II X6 1090T\\
  bogomips: 6421.34\\
  MemTotal: 8191520 kB
\end{quote}

All tests were done on \texttt{openSUSE} 13.1 (x86\_{}64), with Linux kernel 3.11.10-25, \texttt{gcc}~4.9.3 and \texttt{perl}~5.18.1.

All tests were run through \polymake version 2.15-beta3 via the respective interfaces.  This creates some overhead, for instance, due to
data conversion.  While \bb is the only implementation which is actually part of \polymake, this does not constitute a
principal technical advantage over the other convex hull codes tested.  The libraries \texttt{cddlib}~0.94h, \texttt{lrslib} 6.0
and \texttt{libnormaliz}~2.99.4 (which contains the same code as version~3.0) are shipped with \polymake under the GNU General Public License (GPL).  As far as \ppl is concerned the
\polymake distribution only comes with a bare interface, \textit{i.e.}, without the \ppl code. We used the \ppl library version~1.1.

The external software packages used via a file based interfaces are \ftit version~1.6.6, \azove version~2.0, \latte version~1.7.3 and \porta version~1.4.1-20090921.

The \texttt{GMP} was configured to use the standard memory allocator \texttt{malloc}.  Employing a different memory
allocator, such as \texttt{TCMalloc}~\cite{tcmalloc}, can have a great impact, in particular, in a multi-threaded setting.  However, the
precise behavior depends on numerous factors.  For instance, we observed that \texttt{TCMalloc} was clearly superior to
\texttt{malloc} in an \texttt{openSUSE}~12.2 environment, while the difference is only marginal in our setup
with \texttt{openSUSE}~13.1, on the same hardware.

The histogram was created with \texttt{MATLAB} \cite{matlab}.

\clearpage
%\begin{landscape}
\section{Tables with detailed computational results}
\strut
\begin{table}[htp]
  \centering
  \caption{Timings (in seconds) for convex hull computations of non-symmetric cut polytopes $\Cut(G_k)$, see Section~\ref{sec:non-symmetric-cut}.}
  \label{tab:non_sym_cut}
  {\small
  \begin{tabular}[t]{r@{\hspace{.3cm}}rrr@{\hspace{.3cm}}d{3}d{3}d{3}d{3}d{3}d{3}d{3}}
        \toprule
      $k$ & $d$ & $n$ & $m$  & \multicolumn{1}{c}{\bb} & \multicolumn{1}{c}{\cdd}& \multicolumn{1}{c}{\lrs} & %
      \multicolumn{2}{c}{\normaliz} & \multicolumn{1}{c}{\porta} & \multicolumn{1}{c}{\ppl}\\
      &     &     &      &   &  &  & \multicolumn{1}{c}{\footnotesize 1 thread} & \multicolumn{1}{c}{\footnotesize 6 threads} & & \\
        \midrule
      0 &  6 &      32 & 20 &  0.00  &   0.005  &    0.005   &    0.001  &    0.001  &    0.065    &    0.001  \\
      1 &  7 &      64 & 22 &  0.03  &   0.015  &    0.060   &    0.001  &    0.002  &    0.089    &    0.002  \\
      2 &  8 &     128 & 24 &  0.16  &   0.050  &    0.937   &    0.003  &    0.005  &    0.107    &    0.004  \\
      3 &  9 &     256 & 26 &  1.25  &   0.184  &   18.195   &    0.006  &    0.010  &    0.246    &    0.008  \\
      4 & 10 &     512 & 28 & 14.41  &   0.700  &  329.093   &    0.014  &    0.021  &    0.916    &    0.018  \\
      5 & 11 &    1024 & 30 & \oom   &   2.877  & 7699.132^* &    0.038  &    0.047  &    7.980    &    0.043  \\
      6 & 12 &    2048 & 32 & \nc    &  12.086  & \nc        &    0.114  &    0.114  &   98.987    &    0.110  \\
      7 & 13 &    4096 & 34 & \nc    &  50.064  & \nc        &    0.390  &    0.293  & 1810.960    &    0.298  \\
      8 & 14 &    8192 & 36 & \nc    & 210.074  & \nc        &    1.425  &    0.841  & 36729.820^* &    0.890  \\
      9 & 15 &   16384 & 38 & \nc    & 974.681  & \nc        &    5.514  &    2.799  & \nc         &    3.182  \\
     10 & 16 &   32768 & 40 & \nc    & \nc      & \nc        &   22.310  &   10.233  & \nc         &   12.680  \\
     11 & 17 &   65536 & 42 & \nc    & \nc      & \nc        &   93.935  &   42.461  & \nc         &   53.547  \\
     12 & 18 &  131072 & 44 & \nc    & \nc      & \nc        &  399.652  &  181.790  & \nc         &  216.302  \\
     13 & 19 &  262144 & 46 & \nc    & \nc      & \nc        & 1681.473  &  952.263  & \nc         & 1148.847  \\
     14 & 20 &  524288 & 48 & \nc    & \nc      & \nc        & 7173.680^*& 4843.120^*& \nc         & 3641.770  \\
      \bottomrule
  \end{tabular}
  }
\end{table}
%\end{landscape}

{\small
%\begin{table}[htp]
  \begin{longtable}[t]{r@{\hspace{1cm}}d{3}d{3}d{3}d{3}d{3}d{3}d{3}}
  \caption{Running times (in seconds) for integer hull computations of Fibonacci knapsack polytopes $F_d(b)$, see Section~\ref{sec:knapsack:hull}.\label{tab:knapsack:timings}}\\
  \endfirsthead
  \caption[]{(continued)}\\
  \endhead
    \multicolumn{8}{c}{\bb}\\
    \toprule
   $d \setminus b$   & \multicolumn{1}{c}{40} & \multicolumn{1}{c}{50} & \multicolumn{1}{c}{60} & \multicolumn{1}{c}{70} & \multicolumn{1}{c}{80} & \multicolumn{1}{c}{90} & \multicolumn{1}{c}{100}\\
   \midrule
    4 & 0.128 & 0.339 &  0.873 &  1.921 &   3.127 &    7.507 &   16.265 \\
    5 & 0.358 & 1.268 &  4.273 & 15.039 &  42.077 &  136.823 &  355.191 \\
    6 & 0.843 & 3.878 & 21.148 & 76.548 & 305.121 & 1071.274 & 3294.703 \\
    \bottomrule\\

    \multicolumn{7}{c}{\bb (random permutation of input)}\\
    \toprule
   $d \setminus b$   & \multicolumn{1}{c}{40} & \multicolumn{1}{c}{50} & \multicolumn{1}{c}{60} & \multicolumn{1}{c}{70} & \multicolumn{1}{c}{80} & \multicolumn{1}{c}{90} & \multicolumn{1}{c}{100}\\
    \midrule
    4 & 0.069 & 0.138 & 0.295 &  0.572 &  1.035 &  2.023 &   3.738 \\
    5 & 0.220 & 0.555 & 1.203 &  2.901 &  6.861 & 16.215 &  35.602 \\
    6 & 0.596 & 1.844 & 4.495 & 11.417 & 27.508 & 59.907 & 152.199 \\
    \bottomrule\\

\pagebreak

    \multicolumn{7}{c}{\bb (vertices first)}\\
    \toprule
   $d \setminus b$   & \multicolumn{1}{c}{40} & \multicolumn{1}{c}{50} & \multicolumn{1}{c}{60} & \multicolumn{1}{c}{70} & \multicolumn{1}{c}{80} & \multicolumn{1}{c}{90} & \multicolumn{1}{c}{100}\\
    \midrule
    4 & 0.026 & 0.052 & 0.107 & 0.239 &  0.448 &  0.933 &  2.080 \\
    5 & 0.054 & 0.140 & 0.363 & 1.188 &  3.457 &  9.689 & 23.175 \\
    6 & 0.109 & 0.253 & 0.813 & 3.193 & 10.093 & 29.070 & 86.741 \\
    \bottomrule\\

    \multicolumn{7}{c}{\cdd}\\
    \toprule
   $d \setminus b$   & \multicolumn{1}{c}{40} & \multicolumn{1}{c}{50} & \multicolumn{1}{c}{60} & \multicolumn{1}{c}{70} & \multicolumn{1}{c}{80} & \multicolumn{1}{c}{90} & \multicolumn{1}{c}{100}\\
    \midrule
    4 & 2.022 &  8.992 &  30.877 &   90.281 &  231.238 &  537.414 & 1173.495   \\
    5 & 5.667 & 30.601 & 126.978 &  447.216 & 1352.025 & 3594.319 & 8737.240^* \\
    6 & 9.392 & 58.682 & 274.606 & 1092.438 & 3787.679 &    -     &    -       \\
    \bottomrule\\

    \multicolumn{7}{c}{\cdd (random permutation of input)}\\
    \toprule
   $d \setminus b$   & \multicolumn{1}{c}{40} & \multicolumn{1}{c}{50} & \multicolumn{1}{c}{60} & \multicolumn{1}{c}{70} & \multicolumn{1}{c}{80} & \multicolumn{1}{c}{90} & \multicolumn{1}{c}{100}\\
    \midrule
    4 &  2.168 &  9.794 &  34.103 &  100.972 &  266.833 &  717.830 & 1600.367 \\
    5 &  6.195 & 33.444 & 139.920 &  555.176 & 1837.662 & 4928.060 &    -     \\
    6 & 10.379 & 64.724 & 313.493 & 1501.114 & 5145.730 &    -     &    -     \\
    \bottomrule\\

    \multicolumn{7}{c}{\cdd (vertices first)}\\
    \toprule
   $d \setminus b$   & \multicolumn{1}{c}{40} & \multicolumn{1}{c}{50} & \multicolumn{1}{c}{60} & \multicolumn{1}{c}{70} & \multicolumn{1}{c}{80} & \multicolumn{1}{c}{90} & \multicolumn{1}{c}{100}\\
    \midrule
    4 &  2.229 &  9.800 &  33.682 &   99.496 &  282.405 &  717.921 & 1482.070 \\
    5 &  6.162 & 33.009 & 137.282 &  579.661 & 1707.650 & 4614.850 &    -     \\
    6 & 10.132 & 63.107 & 307.064 & 1448.336 & 4777.220 &    -     &    -     \\
    \bottomrule\\

    \multicolumn{7}{c}{\lrs}\\
    \toprule
   $d \setminus b$   & \multicolumn{1}{c}{40} & \multicolumn{1}{c}{50} & \multicolumn{1}{c}{60} & \multicolumn{1}{c}{70} & \multicolumn{1}{c}{80} & \multicolumn{1}{c}{90} & \multicolumn{1}{c}{100}\\
    \midrule
    4 &  5.949 &  32.270 &  132.774   &  448.725   & 1262.392 & 3244.181 & 7583.130^* \\
    5 & 27.849 & 213.319 & 1163.065   & 4957.700^* &    -     &    -     &    -       \\
    6 & 63.492 & 648.386 & 4527.030^* &    -       &    -     &    -     &    -       \\
    \bottomrule\\

    \multicolumn{7}{c}{\lrs (random permutation of input)}\\
    \toprule
   $d \setminus b$   & \multicolumn{1}{c}{40} & \multicolumn{1}{c}{50} & \multicolumn{1}{c}{60} & \multicolumn{1}{c}{70} & \multicolumn{1}{c}{80} & \multicolumn{1}{c}{90} & \multicolumn{1}{c}{100}\\
    \midrule
    4 & 0.385 &  1.317 &   3.829 &   8.170 &   20.089 &   39.051   &  78.252 \\	
    5 & 2.108 &  9.002 &  32.371 & 111.586 &  287.381 &  599.524   & 987.394 \\	
    6 & 5.343 & 32.730 & 170.878 & 644.030 & 1624.175 & 4069.540^* &   -     \\
    \bottomrule\\

  \pagebreak

    \multicolumn{7}{c}{\lrs (vertices first)}\\
    \toprule
   $d \setminus b$   & \multicolumn{1}{c}{40} & \multicolumn{1}{c}{50} & \multicolumn{1}{c}{60} & \multicolumn{1}{c}{70} & \multicolumn{1}{c}{80} & \multicolumn{1}{c}{90} & \multicolumn{1}{c}{100}\\
    \midrule
    4 &  6.765 &  40.089 &  160.902   &  549.346   & 1837.019 & 5302.860^* &  -  \\
    5 & 32.923 & 259.961 & 1516.097   & 7975.330^* &    -     &    -       &  -  \\
    6 & 75.724 & 818.549 & 6345.470^* &    -       &    -     &    -       &  -  \\
    \bottomrule\\

    \multicolumn{7}{c}{\nmz}\\
    \toprule
   $d \setminus b$   & \multicolumn{1}{c}{40} & \multicolumn{1}{c}{50} & \multicolumn{1}{c}{60} & \multicolumn{1}{c}{70} & \multicolumn{1}{c}{80} & \multicolumn{1}{c}{90} & \multicolumn{1}{c}{100}\\
    \midrule

    4 & 0.275 & 0.805 &  2.103 &  5.622 &  13.681 &  29.505 &   65.861 \\
    5 & 0.565 & 2.377 &  8.096 & 26.586 &  74.743 & 194.972 &  421.936 \\
    6 & 0.949 & 3.910 & 17.758 & 68.364 & 184.154 & 438.641 & 1100.974 \\
    \bottomrule\\

    \multicolumn{7}{c}{\porta}\\
    \toprule
   $d \setminus b$   & \multicolumn{1}{c}{40} & \multicolumn{1}{c}{50} & \multicolumn{1}{c}{60} & \multicolumn{1}{c}{70} & \multicolumn{1}{c}{80} & \multicolumn{1}{c}{90} & \multicolumn{1}{c}{100}\\
    \midrule

    4 & 0.070 & 0.180 & 0.410 &  0.965 &  1.734 &   3.495 &   7.696\\
    5 & 0.195 & 0.685 & 1.868 &  5.363 & 12.614 &  37.579 &  96.675\\
    6 & 0.416 & 1.471 & 5.581 & 21.792 & 64.878 & 180.298 & 521.574\\
    \bottomrule\\

    \multicolumn{7}{c}{\porta (random permutation of input)}\\
    \toprule
   $d \setminus b$   & \multicolumn{1}{c}{40} & \multicolumn{1}{c}{50} & \multicolumn{1}{c}{60} & \multicolumn{1}{c}{70} & \multicolumn{1}{c}{80} & \multicolumn{1}{c}{90} & \multicolumn{1}{c}{100}\\
    \midrule

    4 & 0.176 & 0.534 &  1.302 &  2.930 &   5.555 &   8.624 &   15.973\\
    5 & 0.533 & 2.208 &  5.117 & 17.558 &  37.111 &  99.375 &  181.517\\
    6 & 1.163 & 5.903 & 22.705 & 78.768 & 245.711 & 656.937 & 1634.647\\
    \bottomrule\\

    \multicolumn{7}{c}{\porta (vertices first)}\\
    \toprule
   $d \setminus b$   & \multicolumn{1}{c}{40} & \multicolumn{1}{c}{50} & \multicolumn{1}{c}{60} & \multicolumn{1}{c}{70} & \multicolumn{1}{c}{80} & \multicolumn{1}{c}{90} & \multicolumn{1}{c}{100}\\
    \midrule

    4 & 0.035 & 0.073 & 0.133 & 0.258 & 0.337 & 0.501 & 0.940\\
    5 & 0.076 & 0.213 & 0.349 & 0.685 & 0.854 & 2.278 & 4.443\\
    6 & 0.169 & 0.311 & 0.793 & 2.073 & 3.141 & 5.174 & 9.992\\
    \bottomrule\\

    \multicolumn{7}{c}{\ppl}\\
    \toprule
   $d \setminus b$   & \multicolumn{1}{c}{40} & \multicolumn{1}{c}{50} & \multicolumn{1}{c}{60} & \multicolumn{1}{c}{70} & \multicolumn{1}{c}{80} & \multicolumn{1}{c}{90} & \multicolumn{1}{c}{100}\\
    \midrule

    4 & 0.028 & 0.050 & 0.096 & 0.181 &  0.302 &  0.524 &   0.945 \\
    5 & 0.048 & 0.123 & 0.274 & 0.647 &  1.361 &  3.472 &   8.084 \\
    6 & 0.077 & 0.215 & 0.647 & 2.074 &  5.531 & 14.517 &  40.975 \\
    \bottomrule
  \end{longtable}
  \bigskip
%\end{table}
}

\clearpage
\begin{table}[htp]
  \centering
    \caption{Running times (in seconds) for Voronoi diagrams of random point sets, see Section~\ref{sec:voronoi}.}
    \label{tab:voronoi:timings}
  \begin{tabular}[t]{r@{\hspace{1cm}}d{3}d{3}d{3}d{3}d{3}d{3}}
    \multicolumn{7}{c}{\bb}\\
    \toprule
    $d \setminus m$  & \multicolumn{1}{c}{500} & \multicolumn{1}{c}{1000} & \multicolumn{1}{c}{1500} & \multicolumn{1}{c}{2000} & \multicolumn{1}{c}{2500} & \multicolumn{1}{c}{3000}\\
    \midrule
    3 &    0.261 &     0.621  &   1.040 &   1.517 &   2.010 &   2.569 \\
    4 &    1.415 &     3.102  &   4.875 &   6.705 &   8.622 &  10.547 \\
    5 &   11.144 &    25.280  &  40.011 &  55.229 &  69.848 &  85.517 \\
    6 &  100.962 &   255.429  & 425.652 & 608.085 & 800.647 & 995.853 \\
    7 & 1318.488 &  4245.600^*&    -    &    -    &    -    &    -    \\
    \bottomrule\\

    \multicolumn{7}{c}{\cdd}\\
    \toprule
   $d \setminus m$   & \multicolumn{1}{c}{500} & \multicolumn{1}{c}{1000} & \multicolumn{1}{c}{1500} & \multicolumn{1}{c}{2000} & \multicolumn{1}{c}{2500} & \multicolumn{1}{c}{3000}\\
    \midrule
    3 &    6.604 &   28.420 &   67.606 &  124.657 &  200.765 &  295.799 \\
    4 &   53.436 &  246.447 &  599.624 & 1124.073 & 1817.252 & 2712.942 \\
    5 &  326.138 & 1686.490 & 4315.200^* &     -    &     -    &     -    \\
    6 & 2679.876 &     -    &     -    &     -    &     -    &     -    \\
    7 &     -    &     -    &     -    &     -    &     -    &     -    \\
    \bottomrule\\

    \multicolumn{7}{c}{\lrs}\\
    \toprule
    $d \setminus m$ & \multicolumn{1}{c}{500} & \multicolumn{1}{c}{1000} & \multicolumn{1}{c}{1500} & \multicolumn{1}{c}{2000} & \multicolumn{1}{c}{2500} & \multicolumn{1}{c}{3000}\\
    \midrule
    3 &   0.496 &    2.149 &   5.034 &    9.147 &   14.741 &   21.264 \\
    4 &   1.955 &    8.677 &  20.541 &   37.901 &   60.886 &   89.300 \\
    5 &  11.744 &   50.578 & 119.551 &  220.404 &  354.748 &  526.344 \\
    6 &  63.601 &  318.346 & 798.688 & 1530.386 & 2537.572 & 3801.890^*\\
    7 & 404.055 & 2063.702 &    -    &     -    &     -    &     -    \\
    \bottomrule\\

    \multicolumn{7}{c}{\normaliz}\\
    \toprule
    $d \setminus m$ & \multicolumn{1}{c}{500} & \multicolumn{1}{c}{1000} & \multicolumn{1}{c}{1500} & \multicolumn{1}{c}{2000} & \multicolumn{1}{c}{2500} & \multicolumn{1}{c}{3000}\\
    \midrule
    3 &    0.605 &    3.683 &   11.179 &  24.950 &   44.928 &   68.525 \\
    4 &    3.108 &   21.365 &   52.959 & 100.007 &  165.561 &  249.134 \\
    5 &   33.374 &  151.977 &  365.866 & 671.886 & 1085.843 & 1644.610 \\
    6 &  219.029 & 1013.544 & 2436.256 &    -    &    -     &    -     \\
    7 & 1287.456 &    -     &     -    &    -    &    -     &    -     \\
    \bottomrule\\

    \multicolumn{7}{c}{\ppl}\\
    \toprule
    $d \setminus m$ & \multicolumn{1}{c}{500} & \multicolumn{1}{c}{1000} & \multicolumn{1}{c}{1500} & \multicolumn{1}{c}{2000} & \multicolumn{1}{c}{2500} & \multicolumn{1}{c}{3000}\\
    \midrule
    3 &   0.289 &    1.355   &   3.757 &    7.876 &   13.947   &  22.495 \\
    4 &   2.508 &   14.480   &  45.456 &  100.803 &  186.210   & 319.268 \\
    5 &  41.547 &  277.779   & 935.003 & 2255.994 & 4334.870^* &    -    \\
    6 & 917.853 & 6762.370^* &    -    &     -    &     -      &    -    \\
    7 &    -    &     -      &    -    &     -    &     -      &    -    \\
    \bottomrule
  \end{tabular}
\end{table}

{\small
\begin{table}[htp]
\centering
\caption{Timings (in sec.) for counting lattice points in $F_5(b)$; Section~\ref{sec:knapsack:integer}.}
\label{tab:knapsack:lattice:rhs}
\begin{tabular}[t]{r@{\hspace{1cm}}rr@{\hspace{1cm}}d{3}d{3}d{3}d{3}d{3}d{3}d{3}}
      \toprule
    $b$ & $n$ & $\ell$  & \multicolumn{1}{c}{\ftit} & \multicolumn{1}{c}{\bbox}& \multicolumn{2}{c}{\latte} & %
    \multicolumn{2}{c}{\normaliz} & \multicolumn{1}{c}{\proj}\\
        &     &         &                           &                          & \multicolumn{1}{c}{\footnotesize default} & \multicolumn{1}{c}{\footnotesize all primal} & \multicolumn{1}{c}{\footnotesize primal} & \multicolumn{1}{c}{\footnotesize dual} & \\
      \midrule
    40 &  6 &    1366 &    0.500 &   0.247 & 0.012 & 0.028 & 0.048 &  0.007 & 0.011 \\
    50 &  6 &    3137 &    3.331 &   0.501 & 0.014 & 0.028 & 0.024 &  0.014 & 0.023 \\
    60 &  6 &    6509 &   22.171 &   1.193 & 0.013 & 0.027 & 0.212 &  0.039 & 0.039 \\
    70 &  6 &   12182 &  116.258 &   2.420 & 0.014 & 0.027 & 0.069 &  0.080 & 0.063 \\
    80 &  6 &   21245 &  515.834 &   4.928 & 0.013 & 0.027 & 0.121 &  0.160 & 0.100 \\
    90 &  6 &   35025 & 1657.244 &   7.766 & 0.014 & 0.030 & 0.197 &  0.300 & 0.152 \\
   100 &  6 &   55157 & \nc      &  12.787 & 0.013 & 0.027 & 0.329 &  0.555 & 0.216 \\
   110 &  6 &   83616 & \nc      &  20.288 & 0.012 & 0.028 & 0.523 &  0.931 & 0.305 \\
   120 &  6 &  122749 & \nc      &  33.606 & 0.015 & 0.028 & 0.745 &  1.499 & 0.414 \\
   130 &  6 &  175306 & \nc      &  49.113 & 0.013 & 0.028 & 1.173 &  2.342 & 0.556 \\
   140 &  6 &  244473 & \nc      &  64.496 & 0.013 & 0.029 & 1.689 &  3.473 & 0.727 \\
   150 &  6 &  333905 & \nc      &  92.070 & 0.013 & 0.028 & 2.304 &  5.120 & 0.943 \\
   160 &  6 &  447757 & \nc      & 132.536 & 0.014 & 0.029 & 3.147 &  7.150 & 1.229 \\
   170 &  6 &  590715 & \nc      & 177.143 & 0.013 & 0.027 & 4.395 & 10.116 & 1.549 \\
   180 &  6 &  768029 & \nc      & 222.824 & 0.014 & 0.029 & 5.553 & 13.457 & 1.907 \\
   190 &  6 &  985546 & \nc      & 291.333 & 0.012 & 0.028 & 7.214 & 18.045 & 2.352 \\
   200 &  6 & 1249741 & \nc      & 387.886 & 0.014 & 0.028 & 9.235 & 24.116 & 2.880 \\
    \bottomrule
  \end{tabular}
  \medskip
\end{table}

\begin{table}[htp]
\centering
\caption{Timings (in sec.) for counting lattice points in $F_d(60)$; Section~\ref{sec:knapsack:integer}.}
\label{tab:knapsack:lattice:d}
\begin{tabular}[t]{r@{\hspace{1cm}}rr@{\hspace{1cm}}d{3}d{3}d{3}d{3}d{3}d{3}d{3}}
      \toprule
    $d$ & $n$ & $\ell$  & \multicolumn{1}{c}{\ftit} & \multicolumn{1}{c}{\bbox}& \multicolumn{2}{c}{\latte} & %
    \multicolumn{2}{c}{\normaliz} & \multicolumn{1}{c}{\proj}\\
        &     &         &                           &                          & \multicolumn{1}{c}{\footnotesize default} & \multicolumn{1}{c}{\footnotesize all primal} & \multicolumn{1}{c}{\footnotesize primal} & \multicolumn{1}{c}{\footnotesize dual} & \\
      \midrule
    4 &  5 &  4008 &   7.722 &  0.236 &    0.007 &     0.009   &  0.011 & 0.019 & 0.011 \\
    5 &  6 &  6509 &  22.173 &  1.209 &    0.015 &     0.027   &  0.193 & 0.036 & 0.040 \\
    6 &  7 &  7853 &  33.795 &  3.861 &    0.030 &     0.113   &  0.070 & 0.047 & 0.090 \\
    7 &  8 &  8165 &  43.212 &  8.397 &    0.077 &     0.504   &  0.148 & 0.054 & 0.153 \\
    8 &  9 &  8171 &  44.388 & 18.088 &    0.188 &     2.120   &  0.287 & 0.054 & 0.233 \\
    9 & 10 &  8171 &  43.496 & 19.418 &    0.507 &    11.207   &  0.365 & 0.060 & 0.319 \\
   10 & 11 &  8171 &  43.577 & 20.616 &    1.315 &    44.764   &  0.434 & 0.064 & 0.413 \\
   11 & 12 &  8171 &  44.137 & 21.919 &    3.521 &   255.896   &  0.502 & 0.066 & 0.524 \\
   12 & 13 &  8171 &  44.183 & 23.121 &    9.238 &  2357.660   &  0.587 & 0.069 & 0.639 \\
   13 & 14 &  8171 &  44.143 & 24.264 &   24.864 & 36124.990^* &  0.689 & 0.073 & 0.767 \\
   14 & 15 &  8171 &  44.275 & 25.290 &   66.012 &     \nc     & 16.999 & 0.074 & 0.908 \\
   15 & 16 &  8171 &  44.953 & 26.550 &  177.424 &     \nc     & 18.982 & 0.077 & 1.053 \\
   16 & 17 &  8171 &  44.967 & 27.808 &  467.912 &     \nc     & 21.059 & 0.082 & 1.210 \\
   17 & 18 &  8171 &  45.149 & 28.907 & 1246.255 &     \nc     & 22.979 & 0.084 & 1.383 \\
   18 & 19 &  8171 &  45.115 & 30.089 & 3282.651 &     \nc     & 25.063 & 0.084 & 1.571 \\
   19 & 20 &  8171 &  45.233 & 31.594 & \nc      &     \nc     & 27.296 & 0.090 & 1.770 \\
   20 & 21 &  8171 &  45.396 & 32.797 & \nc      &     \nc     & 29.629 & 0.092 & 1.975 \\
    \bottomrule
  \end{tabular}
\end{table}
}

\clearpage

\begin{table}[hptb]
  \centering
  \caption{Running times (in sec) for counting lattice points in random polytopes; Section~\ref{sec:rbox}.}
  \label{tab:randombox:integertimings:}
  \begin{tabular}{r@{\hspace{1cm}}d{3}d{3}d{3}d{3}d{3}d{3}d{3}}
    \multicolumn{7}{c}{\bbox}\\
    \toprule
   $d \setminus b$   & \multicolumn{1}{c}{20} & \multicolumn{1}{c}{30} & \multicolumn{1}{c}{40} & \multicolumn{1}{c}{50} & \multicolumn{1}{c}{60} & \multicolumn{1}{c}{70}\\
    \midrule
    4 &   0.082 &    0.130 &    0.154 &    0.190 &    0.207 &   0.213 \\
    5 &   0.756 &    1.439 &    2.548 &    3.116 &    3.766 &   4.845 \\
    6 &   4.298 &   11.853 &   18.827 &   34.396 &   47.261 &  63.581 \\
    7 &  20.466 &   52.392 &  144.850 &  236.414 &  400.435 & 581.584 \\
    8 &  89.321 &  270.244 &  679.447 & 1693.210 & 2929.315 & -       \\
    9 & 473.808 & 1005.608 & 3511.637 & -        & -        & -       \\
    \bottomrule\\

    \multicolumn{7}{c}{\latte}\\
    \toprule
   $d \setminus b$   & \multicolumn{1}{c}{20} & \multicolumn{1}{c}{30} & \multicolumn{1}{c}{40} & \multicolumn{1}{c}{50} & \multicolumn{1}{c}{60} & \multicolumn{1}{c}{70}\\
    \midrule
    4 &   1.798 &   1.639 &   1.419 &   1.357 &   1.405 &   1.172 \\
    5 & 438.684 & 612.918 & 534.548 & 588.029 & 507.007 & 404.041 \\
    6 &    -    &    -    &   -     &   -     &   -     &   -     \\
    7 &    -    &    -    &   -     &   -     &   -     &   -     \\
    8 &    -    &    -    &   -     &   -     &   -     &   -     \\
    9 &    -    &    -    &   -     &   -     &   -     &   -     \\
    \bottomrule\\

    \multicolumn{7}{c}{\latte (all primal)}\\
    \toprule
   $d \setminus b$   & \multicolumn{1}{c}{20} & \multicolumn{1}{c}{30} & \multicolumn{1}{c}{40} & \multicolumn{1}{c}{50} & \multicolumn{1}{c}{60} & \multicolumn{1}{c}{70}\\
    \midrule
    4 &     0.252   &    0.375 &    0.490   &   0.512 &   0.569 &   0.672 \\
    5 &     3.111   &    4.651 &    6.690   &   8.114 &   9.516 &  11.017 \\
    6 &    49.751   &  122.619 &  154.804   & 235.039 & 281.204 & 338.060 \\
    7 &   801.512   & 3016.324 & 4459.500^* &    -    &    -    &    -    \\
    8 & 22332.390^* &     -    &     -      &    -    &    -    &    -    \\
    9 &      -      &     -    &     -      &    -    &    -    &    -    \\
    \bottomrule\\

    \multicolumn{7}{c}{\normaliz}\\
    \toprule
   $d \setminus b$   & \multicolumn{1}{c}{20} & \multicolumn{1}{c}{30} & \multicolumn{1}{c}{40} & \multicolumn{1}{c}{50} & \multicolumn{1}{c}{60} & \multicolumn{1}{c}{70}\\
    \midrule
    4 & 0.002 &  0.003 &   0.004 &   0.006 &   0.006 &   0.005 \\
    5 & 0.005 &  0.007 &   0.013 &   0.020 &   0.022 &   0.025 \\
    6 & 0.017 &  0.040 &   0.085 &   0.129 &   0.175 &   0.237 \\
    7 & 0.098 &  0.426 &   0.782 &   1.307 &   1.941 &   2.392 \\
    8 & 1.068 &  3.214 &   8.565 &  16.734 &  24.145 &  34.041 \\
    9 & 5.256 & 50.533 & 129.305 & 234.613 & 474.163 & 626.180 \\
    \bottomrule\\

%    \multicolumn{7}{c}{\normaliz (dual mode)}\\
%    \toprule
%   $d \setminus b$   & \multicolumn{1}{c}{20} & \multicolumn{1}{c}{30} & \multicolumn{1}{c}{40} & \multicolumn{1}{c}{50} & \multicolumn{1}{c}{60} & \multicolumn{1}{c}{70}\\
%    \midrule
%    4 &   0.103 &    0.200 & 0.171 & 0.221 & 0.220 & 0.231 \\
%    5 & 110.524 & 2314.403 &  -    &  -    &  -    &  -    \\
%    6 &    -    &     -    &  -    &  -    &  -    &  -    \\
%    7 &    -    &     -    &  -    &  -    &  -    &  -    \\
%    8 &    -    &     -    &  -    &  -    &  -    &  -    \\
%    9 &    -    &     -    &  -    &  -    &  -    &  -    \\
%    \bottomrule\\

    \multicolumn{7}{c}{\projection}\\
    \toprule
   $d \setminus b$   & \multicolumn{1}{c}{20} & \multicolumn{1}{c}{30} & \multicolumn{1}{c}{40} & \multicolumn{1}{c}{50} & \multicolumn{1}{c}{60} & \multicolumn{1}{c}{70}\\
    \midrule
    4 &    0.013 &    0.023 &    0.029 &    0.035 &    0.042 &    0.040 \\
    5 &    0.100 &    0.232 &    0.462 &    0.621 &    0.762 &    1.032 \\
    6 &    0.450 &    1.806 &    3.945 &    6.306 &   10.476 &   14.247 \\
    7 &    1.369 &    9.502 &   28.062 &   57.688 &  102.579 &  149.536 \\
    8 &    3.414 &   41.705 &  158.855 &  138.180 &  148.227 &  223.361 \\
    9 &    7.490 &  151.081 &  219.040 &  248.157 &  769.460 & 1017.670 \\
    \bottomrule
  \end{tabular}
\end{table}

  \bigskip

{
\begin{table}[htp]
\centering
\caption{Timings (in seconds) for counting lattice points in the matching polytope of $K_n$, see Section~\ref{sec:matching}.}
\label{tab:matching:lattice}
\begin{tabular}[t]{rr@{\hspace{.5em}}d{3}d{3}d{3}d{3}d{3}d{3}d{3}d{3}}
      \toprule
    $n$ & $m$ & \multicolumn{1}{c}{\ftit} & \multicolumn{1}{c}{\azove} & \multicolumn{1}{c}{\bbox}& \multicolumn{2}{c}{\latte} & %
    \multicolumn{2}{c}{\normaliz} & \multicolumn{1}{c}{\proj}\\
       &     &         &                           &                          & \multicolumn{1}{c}{\footnotesize default} & \multicolumn{1}{c}{\footnotesize all primal} & \multicolumn{1}{c}{\footnotesize primal} & \multicolumn{1}{c}{\footnotesize dual} & \\
      \midrule
    4 &  10  &    0.00 &     0.000   &    0.005   &   0.028 &  0.070 &    0.000 &   0.001 &   0.002\\
    5 &  15  &    0.00 &     0.000   &    0.051   &   0.450 & 16.865 &    0.006 &   0.003 &   0.009\\
    6 &  21  &    0.00 &     0.003   &    2.430   &  13.195 &  \nc   &    0.093 &   0.008 &   0.059\\
    7 &  28  &    0.18 &     0.004   &  281.244   & 531.014 &  \nc   & 1782.251 &   0.072 &   0.507\\
    8 &  36  &   12.05 &     0.004   & 6138.776^* &   \nc   &  \nc   &  \nc     &   2.163 &   4.306\\
    9 &  45  & 1557.38 &     0.004   &  \nc       &   \nc   &  \nc   &  \nc     &  95.743 &  60.759\\
   10 &  55  & \nc     &     0.009   &  \nc       &   \nc   &  \nc   &  \nc     & 558.374 &  \nc   \\
   11 &  66  & \nc     &     0.017   &  \nc       &   \nc   &  \nc   &  \nc     &   \nc   &  \nc   \\
   12 &  78  & \nc     &     0.062   &  \nc       &   \nc   &  \nc   &  \nc     &   \nc   &  \nc   \\
   13 &  91  & \nc     &     0.236   &  \nc       &   \nc   &  \nc   &  \nc     &   \nc   &  \nc   \\
   14 & 105  & \nc     &     1.046   &  \nc       &   \nc   &  \nc   &  \nc     &   \nc   &  \nc   \\
   15 & 120  & \nc     &     4.696   &  \nc       &   \nc   &  \nc   &  \nc     &   \nc   &  \nc   \\
   16 & 136  & \nc     &    23.224   &  \nc       &   \nc   &  \nc   &  \nc     &   \nc   &  \nc   \\
   17 & 153  & \nc     &   112.785   &  \nc       &   \nc   &  \nc   &  \nc     &   \nc   &  \nc   \\
   18 & 171  & \nc     &   578.191   &  \nc       &   \nc   &  \nc   &  \nc     &   \nc   &  \nc   \\
   19 & 190  & \nc     &  2916.652   &  \nc       &   \nc   &  \nc   &  \nc     &   \nc   &  \nc   \\
   20 & 210  & \nc     & 15462.050^* &  \nc       &   \nc   &  \nc   &  \nc     &   \nc   &  \nc   \\
    \bottomrule
  \end{tabular}
\end{table}
}

\FloatBarrier

%\clearpage
\bibliographystyle{amsplain}
\bibliography{MPC}
\clearpage

\end{document}